%
%
%
%
%
\input amstex
\documentstyle{amsppt}
\magnification =1200
\pagewidth {6.5 true in}
\pageheight {8.5 true in }
\hcorrection {-.2 true in}
\vcorrection {+.1 true in}

\baselineskip = 14pt
\parindent 20 pt
\parskip = 3pt plus 1pt minus .5pt


\define\ds{\displaystyle}
\define\fp{\flushpar}
\define\cstar{$\roman{C}^*$}   
   
\font\bigmath=cmsy10 scaled \magstep1
\def\bigtimes{\mathbin{\hbox{\bigmath\char'2}}}

\define\presen#1#2#3#4#5{{#1}_{#3} @>>> {#1}_{#4} @>>>
   {#1}_{#5} @>>> \cdots @>>> #2}  

\define\Tpresen#1#2#3#4{T_{#1} @>>> T_{#2} @>>> T_{#3} 
   @>>> \cdots @>>> {#4}} 

\define\algtri{(D, A, B)}
\define\spectri{(X, P, G)}
\define\BB{\Cal B} 

\define\II{\Cal I}
\define\JJ{\Cal J}
\define\KK{\Cal K}
\define\MM{\Cal M}
\define\UU{\Cal U}
\define\orb{\text{orb}}
\define\s{\sigma}
\redefine\t{\tau}
\define\sab{\s_{a,b}}
\define\tab{\tau_{a,b}}
\define\Id{\roman{Id}}
\define\sc{\overline{P \cap (I \bigtimes I)}}
\define\tsub{\text{sub}}
\define\hs{\ell^2(I)}
\define\ehat{{\hat e}}
\define\fhat{{\hat f}}

\define\xhat{{\hat x}}
\define\phst{^{\phantom{*}}}
\define\dx{\delta_x}
\define\dy{\delta_y}
\define\ear{\tilde{R}_{\infty}(c)}
\define\MIC{\roman{MIC}}
\define\dist{\operatorname{dist}}

\topmatter
\title Meet Irreducible Ideals in Direct Limit Algebras \endtitle
\author
Allan P. Donsig \\
Alan Hopenwasser \\
Timothy D. Hudson$^*$ \\
Michael P. Lamoureux$^\dag$ \\
Baruch Solel 
\endauthor
\leftheadtext{Donsig, Hopenwasser, Hudson, Lamoureux, \& Solel}
\address
Pure Mathematics Dept., University of Waterloo, Waterloo, ON,
Canada, N2L 3G1
\endaddress
\email
apdonsig\@math.uwaterloo.ca
\endemail
\address
Dept. of Mathematics, University of Alabama, Tuscaloosa,
AL, U.S.A., 35487
\endaddress
\email
ahopenwa\@gp.as.ua.edu
\endemail
\address
Dept. of Mathematics, East Carolina University, Greenville,
NC, U.S.A., 27858 
\endaddress
\thanks $^*$Partially supported by an NSF grant. \endthanks
\email
tdh\@math1.math.ecu.edu
\endemail
\address 
Dept. of Mathematics and Statistics, University of Calgary, 
2500 University Drive N.W., Calgary, AB, Canada T2N 1N4
\endaddress
\email
mikel\@math.ucalgary.ca
\endemail
\thanks $^\dag$Partially supported by an NSERC grant. \endthanks
\address
Dept. of Mathematics, Technion--Israel Institute of Technology,
32000 Haifa, Israel
\endaddress
\email
mabaruch\@tx.technion.ac.il
\endemail
\endtopmatter
   
\document

We study the meet irreducible ideals (ideals $I$ so that $I = J \cap K$ 
implies $I=J$ or $I=K$) in certain direct limit algebras.
The direct limit algebras will generally be strongly maximal triangular 
subalgebras of AF \cstar-algebras, or briefly, strongly maximal TAF algebras.
Of course, all ideals are closed and two-sided.

These ideals have a description in terms of the coordinates, or spectrum, 
that is a natural extension of one description of meet irreducible ideals 
in the upper triangular matrices.
Additional information is available if the limit algebra is an analytic 
subalgebra of its \cstar-envelope or if the analytic algebra is trivially 
analytic with an injective 0-cocycle.
In the latter case, we obtain a complete description of the
meet irreducible ideals, modeled on the description
in the algebra of upper triangular matrices.  
This applies, in particular, to all full nest algebras.

One reason for interest in the meet irreducible ideals of a strongly maximal
TAF algebra is that each meet irreducible ideal is the kernel of 
a nest representation of the algebra (Theorem 2.4).
A {\it nest representation\/} of an operator algebra $A$ is a 
norm continuous representation of $A$ acting on a Hilbert space 
with the property that the lattice of closed invariant subspaces 
for the representation is totally ordered. 
These representations were introduced in \cite{L1} as analogues
for a general operator algebra of 
the irreducible representations of a \cstar-algebra.
The meet irreducible ideals seem analogous to the primitive ideals in 
a \cstar-algebra.
Indeed, in a \cstar-algebra, the meet irreducible ideals are precisely
the primitive ideals \cite{L3, Theorem 2.1}.

This analogy can be extended by noting that the meet irreducible ideals 
form a topological space under the hull-kernel topology and every ideal 
is the intersection of the meet irreducible ideals which contain it.  
There is a one-to-one correspondence between closed sets in the meet 
irreducible ideal space and ideals in the strongly maximal TAF algebra; 
thus the topological space of meet irreducible ideals determines 
completely the lattice of ideals of the limit algebra,
exactly as the primitive ideal space does for \cstar-algebras.
Similar results have been obtained for other operator algebras, 
including the compacts in a nest algebra, the disc algebra, and various
nonselfadjoint crossed products \cite{L1,L2,L3}.

An interesting subset of the meet irreducible ideals are the completely 
meet irreducible ideals, namely those satisfying an analogous condition, 
only for arbitrary intersections instead of just for finite intersections.
We describe these ideals and show that, for direct limit algebras 
generated by their order preserving normalizers, this subset is isomorphic
to the spectrum of the limit algebra (Theorem 5.3).
Also, there is a distance formula for ideals in a strongly maximal TAF
algebra (Theorem 6.2) that is analogous to Arveson's distance formula 
for nest algebras and to the distance formulae in \cite{MS2}.

\head  0. Algebras \& Coordinates  \endhead

An analysis of ideals in direct limit algebras is greatly
facilitated by the technique of coordinatization.
After outlining the algebraic setting, we describe the essential 
ingredients for coordinatization in the context in which we need it; 
for more detail on coordinatizations and more general results 
the reader is referred to \cite{R}, \cite{MS1}, and \cite{P4}.
The book \cite{P4} by Power is also a convenient reference for 
direct limit algebras.

If $A$ is a strongly maximal triangular subalgebra of a unital
AF \cstar-algebra, $B$, then $D = A \cap A^{*}$ is a 
canonical masa in $B$ and $A + A^{*}$ is dense in $B$.
(This is one definition of ``strongly maximal triangular''.)
Since $B$ is AF, it may be written as a direct limit of finite
dimensional \cstar-algebras:
$$ 
\presen B B 1 2 3.
$$
In turn, $A$ can be written as a direct limit
$$
\presen A A 1 2 3,
$$ 
where each $A_n$ is a maximal triangular subalgebra
of $B_n$; also, $D$ is a direct limit
$$
\presen D D 1 2 3,
$$
where each $D_n = A_n \cap A_n^{*}$ is a masa in $B_n$.
If $I$ is a two-sided ideal of $A$ then, essentially because
$I$ is a $D$-bimodule, it follows that $I$ is the closed union
of the $I \cap A_n$.
\par
Furthermore, it is possible to select a system of matrix
units for $B$ so that each of $A$ and $D$ are generated
by the matrix units which they contain.  
Of course, it follows that every ideal is also generated by
the matrix units it contains.
The system of matrix units can also be chosen so that each 
matrix unit in $B_n$ is a sum of matrix units in $B_{n+1}$.  
By identifying each $B_n$ with its natural image in $B$, 
we may consider all the embeddings which appear in the direct system 
to be inclusions.
\par
A direct system whose limit is $A$ will be referred to as a
{\it presentation\/} for $A$.  
Given a presentation for $A$ as above, we can construct another
presentation by choosing a subsequence $A_{n_1},A_{n_2},\ldots$,
with maps given by composing the maps from the original presentation.
We call this new presentation a {\it contraction} of the original.
\par
We now coordinatize the triple of algebras $\algtri$,
where $B$ is an AF \cstar-algebra and $A$ is a strongly maximal triangular
subalgebra of $B$ whose diagonal is $D$.  
Also assume that a system of matrix units for $B$ has been selected with
the properties described above.  
We need to define a {\it spectral triple\/} $\spectri$ for $\algtri$.  
The first ingredient, $X$, is simple: it is just the maximal ideal
space for $D$.  
So $D$ is isomorphic to $C(X)$ and, since $D$ is a direct limit of 
finite dimensional algebras, $X$ is isomorphic to the Cantor set.
\par
Since the \cstar-algebra, $B$, is AF, it is a groupoid algebra;
$G$ will be the groupoid associated with $B$.  While we will use
some of the language of groupoids and a couple of results about
groupoids, the reader does not need extensive knowledge of 
groupoids in order to follow our arguments.  Indeed, $G$ is a 
special type of groupoid and we can describe it completely in
a very naive fashion.  Each matrix unit, $e$, from the system
of matrix units for $B$ acts on $D$ by conjugation 
($e^*De \subseteq D$); consequently, each matrix unit, $e$,
induces a partial homeomorphism of $X$ into itself
(i.e., a homeomorphism between two clopen subsets of $X$).
Let $\ehat$ denote the graph of this homeomorphism.
\par
As a set, $G$ is simply the union of the graphs of all the
partial homeomorphisms induced by matrix units.
Thus, $G$ is a subset of $X \bigtimes X$; it is not difficult
to check that it is an equivalence relation.  
There is, however, an additional structure, a topology, on $G$.  
This topology is the smallest topology in which every $\ehat$ is an 
open subset.
It turns out that every $\ehat$ is also closed, and hence compact.
This description of $G$ appears to be dependent on the choice
of matrix unit system (and hence on the choice of presentation);
in point of fact the same topological equivalence relation
arises from any choice of presentation and any choice of matrix
unit system.  
Indeed, in place of matrix units one may use the collection 
of all partial isometries in $B$ which normalize $D$.
(A partial isometry, $v$, {\it normalizes} $D$ if $vDv^* \subseteq D$
and $v^*Dv \subseteq D$.)
\par
A topological equivalence relation such as $G$ is an $r$-discrete,
principal, topological groupoid.  We won't use all this terminology,
but we do need to say what the groupoid operations are.  Two
elements $(x,y)$ and $(w,z)$ are composable if, and only if,
$y=w$.  In that case, the product is given by $(x,y) \circ (y,z)
 = (x,z)$.  Inverses are given by $(x,y)^{-1} = (y,x)$.
\par
The graph, $\nu$, of the partial homeomorphism associated with
a matrix unit (or with a normalizing partial isometry) has 
the following properties:
\item{i)} if $(x,y_1) \in \nu$ and $(x, y_2) \in \nu$, then
 $y_1 = y_2$,
\item{ii)} if $(x_1, y) \in \nu$ and $(x_2, y) \in \nu$, then
 $x_1 = x_2$.
\fp
A subset of $G$ with these properties is called a {\it $G$-set\/}.
It is a property of the topology on $G$ that any point has
a neighborhood basis which consists of compact, open $G$-sets.
All $G$-sets which appear in this paper can be taken to to be
compact and open; assume that any $G$-set which appears is compact 
and open even if these adjectives are absent.
\par
If $\nu_1$ and $\nu_2$ are  $G$-sets, then so
is the composition $\nu_1 \circ \nu_2$, which is defined to 
be the set $\{ a\circ b \: a \in \nu_1, b \in \nu_2 
\text{ and } a \text{ and } b \text{ are composable\,}\}$.
In the case of graphs $\ehat$ and $\fhat$ of matrix units
(or normalizing partial isometries), $\ehat \circ \fhat$
will be the graph of the product $ef$ in $B$.
\par
The space, $X$, can be identified with the diagonal of $G$
via the homeomorphism $x \leftrightarrow (x,x)$.  In particular,
the diagonal of $G$ is an open subset of $G$.  (For readers 
familiar with groupoids, the diagonal is the space, $G^0$,
of units of $G$.  The fact that it is open means that $G$ is
$r$-discrete.)  One should also note that in the present 
context, the two coordinate projection maps $\pi_1$ and
$\pi_2$, when restricted to $G$, are open maps (from the
groupoid topology on $G$ to the topology on $X$); in fact, they
are local homeomorphisms.
\par
It remains to describe the middle component, $P$, of the spectral
triple.  The short way is to invoke the spectral theorem for
bimodules \cite{MS1}: $P$ is the unique open subset of $G$ which is
the support set for the subalgebra $A$.  The fact that $A$ is generated
by the matrix units which it contains permits a naive definition
of $P$: it is simply the union of the graphs, $\ehat$, for the
matrix units $e$ in $A$.  As such, it is a subrelation of $G$ and it
carries the relative topology induced by the topology on $G$.
The apparent dependence of $P$ on choice of matrix unit system
(or presentation) is illusory and $P$ is, in fact, an isometric
isomorphism invariant for $A$ \cite{P2}. 
 We shall call $P$ the {\it spectrum\/} of $A$.
\par
As is to be expected, properties of $A$ are reflected in
properties of $P$.  The fact that $A$ is an algebra means
that $P \circ P \subseteq P$.  The triangularity of $A$
becomes the property that $P \cap P^{-1}$ is the diagonal
of $G$.  Finally, strong maximality for $A$ is equivalent
to $P \cup P^{-1} = G$.  Note, in particular, that the
topological relation, $P$, induces a total order on each
equivalence class in $G$.  We shall need a notation for
equivalence classes: if $z \in X$, let $\orb_z = \{
x \in X \: (x,z) \in G\,\}$.  We sometimes emphasize the
induced order on each equivalence class by writing
$x \le y$ when $(x,y) \in P$.
\par
Some of our results are valid in the context of triangular
subalgebras of $B$ which are analytic.  The simplest definition
of analytic subalgebras is in terms of real valued cocycles.
A continuous function $c\: G \longrightarrow \Bbb R$ is a
{\it 1-cocycle} provided that $c(x,y) + c(y,z) = c(x,z)$, for all 
$(x,y), (y,z) \in G$.  We say that $A$ is {\it analytic\/}
if $P = c^{-1}[0, \infty)$.  We say that $A$ is {\it trivially
analytic\/} when $c$ has the special form $c(x,y) = b(y) - b(x)$
for a continuous function $b\:X \longrightarrow \Bbb R$.  (Such
a function, $b$, is called a {\it 0-cocycle} and $c$ is the 
{\it coboundary} of $b$.)  
The material in Section 3 will be valid for trivially
analytic algebras with the additional requirement that the
0-cocycle be an injective function.  This family of algebras includes
all full nest algebras.  (While the 0-cocycle most naturally
associated with a full nest algebra will not be injective, it can
be replaced by an injective 0-cocycle whose coboundary yields the 
same analytic algebra.)
\par
Just as the algebra $A$ has a natural support set $P \subset G$, 
each two sided closed ideal $\II \subseteq A$ has a support set,
$\s$.  The existence of $\s$ is given by the spectral theorem for
bimodules and a complete description of coordinatization for ideals
is given in \cite{MS1}.  Also, just as before, a naive description
of $\s$ is available based on the fact that an ideal is generated
by the matrix units which it contains \cite{P1}.  
So, $\s$ is the union of
the graphs associated with matrix units of $\II$ and the topology
is the relative topology from $P$.  The definition of $\s$ is, of
course, independent of choice of matrix unit system or presentation.
\par
The fact that $\II$ is an ideal is reflected in the following
property for~$\s$: 
$$
(w,x) \in P, (x,y) \in \s, (y,z) \in P \Longrightarrow
 (w,z) \in \s.
$$
We say that an open subset of $P$ which satisfies this property
is an {\it ideal set\/}.  Lemma 4.3 in \cite{MS1} shows that 
each ideal set is the support set of a closed, two sided ideal 
in $A$.
\par
We say that an ideal set, $\s_1$, is {\it meet irreducible\/}
if, whenever $\s = \t_1 \cap \t_2$ with $\t_1, \t_2$ ideal sets,
either $\s = \t_1$ or $\s = \t_2$.  Since intersection of ideals
corresponds to intersection of ideal sets, an ideal in $A$ is 
meet irreducible if, and only if, the corresponding ideal set
is meet irreducible.
\par
\head 1. MI-chains \endhead

In $T_n$, the algebra of $n \times n$ upper triangular matrices,
each meet irreducible ideal is determined by a matrix	
unit.  If $e_{st}$ is a matrix unit in $T_n$, then the	
meet irreducible ideal $\II$
 associated with $e_{st}$ is the
largest ideal in $T_n$ which does not contain $e_{st}$.
This ideal is generated as a linear subspace by	the
set of matrix units $e_{nm}$ where either $n < s$	
or $m > t$.

The meet irreducible ideals in $T_n$ can also be
described in terms of the coordinates, rather 
than  matrix units.  Let $X = \{1, \dots, n\}$ and 
$P = \{(s,t) \: s,t \in X \text{ and } s \le t \}$.  
Then $P$ is the support set for $T_n$.  
Let $I$ be an interval contained in $X$; 
i.e., $I = \{i \: s \le i \le t \}$ for some $s,t$.  
Then the meet irreducible ideal $\II$ associated with $I$ is 
the set of all	matrices supported on 
$P \setminus P \cap (I \bigtimes I)$.

In the TAF algebra context, the description of meet 
irreducible ideals in terms of coordinates needs
almost no modification from the finite dimensional case.  
The description in terms of matrix units is considerably 
more complicated than the finite dimensional description.  
In [La], Lamoureux gave a procedure for constructing meet irreducible 
ideals from certain sequences of matrix units, provided that
the embeddings satisfy a special condition.  
(This condition is met by standard embeddings, by refinement
embeddings, and, more generally, by nest embeddings.)
However, this procedure fails to give all meet 
irreducible ideals even in the simplest TAF algebras.  
\par
There is, in fact, a more general family of matrix unit sequences 
from which  meet irreducible ideals can be constructed.  
This concept -- MI-chains of matrix units -- yields all 
the meet irreducible ideals (provided we consider all possible 
contractions of a given presentation); 
furthermore, it is valid for all TAF algebras.
\par
Let $A$ be a TAF algebra with presentation
$$
\presen A A {1} {2} {3}.
$$
\definition{Notation} If $e \in A_n$, then $\Id_n(e)$
will denote the ideal generated by $e$ in $A_n$.
If $k >n$, then $e \in A_k$ also; therefore
$\Id_k(e)$ is defined and $\Id_n(e) \subseteq \Id_k(e)$.
\enddefinition
\definition{Definition {1.1}}
A sequence $(e_k)_{k \ge N}$ of matrix units from $A$ 
will be called an {\it MI-chain\/} if the following
two conditions are satisfied for all $k \ge N$:
\item{(A)} $\ds e_k \in A_k$.
\item{(B)} $e_{k+1} \in \Id_{k+1}(e_k)$.
\enddefinition 
\par
If $(e_k)$ is an MI-chain for $A$, let $\II$ be the join
of all ideals which do not contain any matrix unit
$e_k$ from the chain.  In other words, $\II$ is the
largest ideal in $A$ which does not contain any $e_k$.
\par
\proclaim{Theorem 1.2} 
Let $A$ be a strongly maximal TAF algebra with some presentation.
For each MI-chain $(e_k)_{k \ge N}$ from the presentation,
the ideal $\II$ associated with $(e_k)$ is meet irreducible.  
Conversely, every proper meet irreducible ideal in $A$ is induced by 
some MI-chain, chosen from some contraction of this presentation.  
\endproclaim
\demo{Proof} Let $(e_k)$ be an MI-chain of matrix units
and let $\II$ be the corresponding ideal.  Suppose that
$\JJ$ and $\KK$ are two ideals each of which properly
contains $\II$.  Since $\II$ is the largest ideal containing
no matrix units from the MI-chain, there exist indices $s$
and $t$ such that $e_s \in \JJ$ and $e_t \in \KK$.
Condition (B) in the definition of MI-chain implies that
$e_n \in \JJ$ for all $n > s$ and $e_m \in \KK$ for
all $m > t$.  Thus, $\JJ \cap \KK$ contains matrix 
units from the MI-chain, which implies that $\JJ
\cap \KK$ properly contains $\II$.  This proves that
$\II$ is meet irreducible.
\par
For the converse, suppose that $\II$ is a proper meet 
irreducible ideal in $A$ and that
$$
\presen A A 1 2 3
$$
is a presentation for $A$.  Each $A_k$ is a maximal triangular
subalgebra of a finite dimensional \cstar-algebra.  Let
$\II_k = \II \cap A_k$, for each k.  While $\II$ is the
closed union of the $\II_k$, it is not necessarily the case 
that each $\II_k$ is meet irreducible as an ideal in $A_k$.
Note that, by contracting the presentation if necessary,
we may also assume that $\II_k$ is a proper ideal in
$A_k$, for each $k$.
\par
 From the known structure of ideals in maximal triangular
subalgebras of finite dimensional \cstar-algebras, it
follows that for each $k$ there is a minimal set $E_k$
of matrix units in $A_k \setminus \II_k$
 such that any ideal of $A_k$
which is larger than $\II_k$ must contain one of the
matrix units in $E_k$.  Begin the construction of an
MI-chain for $\II$ by letting $e_1$ be any matrix
unit from $E_1$. 
\par
For each $e \in E_1$, let $\JJ_e$ denote the ideal in
$A$ generated by $\II$ and $e$.  Since each such $e$
is not in $\II_1$ but is in $A_1$, $e$ does
not belong to $\II$; thus
 $\II$ is a proper
subset of each $\JJ_e$.  Let $\JJ = \cap \{\JJ_e \:
e \in E_1 \}$.  This is a finite intersection and
$\II$ is meet irreducible, so $\JJ$ properly contains $\II$.  
Consequently, for some $k \ge 2$, $\JJ \cap A_k$ 
properly contains $\II_k = \II \cap A_k$.  
By replacing the presentation by a contraction and 
relabeling, we may assume that $k=2$.
\par
Since $\JJ \cap A_2$ properly contains $\II_2$, there is a 
matrix unit $e_2 \in E_2$ such that $e_2 \in \JJ \cap A$.
By the construction of $\JJ$, $e_2 \in \Id_2(e_1)$;  
 thus condition (B) for MI-chains is satisfied by
the pair $e_1$, $e_2$.  
\par
If we now iterate this construction, we obtain a presentation
which is a contraction of the original presentation and a
sequence of matrix units $(e_k)_{k \ge 1}$ which is an
MI-chain.  Since $\II$ contains none of the $e_k$, $\II$
is a subset of the meet irreducible ideal associated with
the MI-chain.  But if $\KK$ is an ideal larger than $\II$, 
then $\KK \cap A_k$ properly contains $\II_k$ for some
$k$ and hence $\KK$ contains some element of $E_k$.
By the construction of the sequence $(e_n)$,
$e_{k+1}$ is in the ideal generated by each element
of $E_k$; hence $e_{k+1} \in \KK$.
Thus, $\II$ is the largest ideal which contains 
none of the $e_k$ and so it is the meet irreducible 
ideal associated with the MI-chain.  \qed
\enddemo
\par
It is natural to ask if there is a 1-1 correspondence between
MI-chains and meet-irreducible ideals.
Without other conditions, the answer is clearly no.
For example, take an MI-chain for the zero ideal and 
change the first finitely many matrix units in the MI-chain.
To fix this trivial kind of counterexample, 
the appropriate condition on the MI-chain is 
\item{(C)} for a matrix unit $f$ in $A_k$, if $f \in \Id_k(e_k)$ 
  and $f \ne e_k$, then $e_{k+1} \notin \Id_{k+1}(f)$.
\fp
In fact, Theorem~1.2 always gives an MI-chain satisfying this condition.
Using the notation of the proof, 
observe that if $f$ is a matrix unit in $A_1$ 
which is not equal to $e_1$ but is in $\Id_1(e_1)$, 
then $f$ belongs to $\II_1$, and hence to $\II$.  
Observe that $\Id_2(f) \subseteq \II$, and so $e_2 \notin \Id_2(f)$.
This verifies condition (C) for the pair $e_1$, $e_2$ and,
by induction, the MI-chain $(e_k)$ satisfies the condition.

However, restricting to MI-chains satisfying condition (C) still 
does not give a 1-1 correspondence, as the following example shows. 
Thus the correspondence between meet irreducible ideals and 
MI-chains is rather subtle.

\definition{Example 1.3}
For $n \ge 1$, let $A_n = T_{2^n} \oplus T_{2^n}$ and 
let $\alpha_n \colon A_n\rightarrow A_{n+1}$ be given by 
the block matrix map
$$ 
\bmatrix A & B \\ & C \endbmatrix \oplus \bmatrix D & E \\ & F \endbmatrix
     \longrightarrow 
     \bmatrix A & & & B \\ & D & E \\ & & F \\ & & & C \\ \endbmatrix
     \oplus
     \bmatrix D & & & E \\ & A & B \\ & & C \\ & & & F \\ \endbmatrix .
$$
Consider the algebra $A$ which is the direct limit of the 
algebras $A_n$ with respect to the maps $\alpha_n$.
For each $n$, let $e_n$ be the upper-right matrix unit of $B$ in each 
$A_n$, and let $f_n$ be the upper-right matrix unit of $E$ in each $A_n$.
Observe that for each $n$, $e_{n+1}$ is a summand of $e_n$, and 
so $e_{n+1}\in\Id_{n+1}(e_n)$. 
Similarly, $f_{n+1}\in\Id_{n+1}(f_n)$, and so both $(e_n)_n$ and 
$(f_n)_n$ are MI-chains.
Moreover, since there is no matrix unit $f$ in $A_k$ with 
$f\in \Id_k(e_k)$ and $f \ne e_k$, then $(e_n)$ satisfies 
condition (C), and similarly for $(f_n)$.
It is easy to see that both chains correspond to the zero ideal 
of $A$. \qed
\enddefinition
\par
\head 2. Meet Irreducible Ideals and Nest Representations  \endhead
In this section we will construct meet irreducible ideals
using coordinate methods.  
Fix notation as follows: 

\noindent
{\bf Notation.} 
Let $A$ be a strongly maximal TAF algebra whose enveloping
\cstar-algebra is $B$ and whose diagonal is $D$, a canonical
masa in $B$.  
Also, $(X, P, G)$ will denote the spectral
triple for $(D, A, B)$.  

For subsets of $G$, the closure operator will always denote 
closure with respect to the groupoid topology on $G$, never 
the relative product topology on the larger set $X \bigtimes X$.
Also, by an order interval in an equivalence class of $G$ we mean
the set of points $\{ y \in X : (x,y),(y,z) \in P \}$, where
$(x,z) \in P$, possibly excluding the endpoints $x$ and $z$.

\proclaim{Theorem 2.1} 
With notation as above, let $I$ be an order interval from 
an equivalence class from $G$ and let $\s = P \setminus \sc$.  
Then $\s$ is a meet irreducible ideal set.
\endproclaim
\demo{Proof} We will first show that $\s$ is an ideal
set in $P$.  
To that end, assume that $(u,x) \in P$ and $(x,y) \in \s$.  
We will show that $(u,y) \in \s$.
\par
Suppose, to the contrary, that $(u,y) \in \sc$.  Then
there are sequences $u_n$ and $y_n$ in $I$ such that
$(u_n, y_n) \in P$ and $(u_n, y_n) \longrightarrow
(u,y)$ in $P$.  Let $T$ and $S$ be compact, open $G$-sets
containing $(u,x)$ and $(x,y)$ respectively. 
We may select $T$ and $S$ so that each is a subset of $P$.
 (These sets 
may be chosen to be the graphs of matrix units in $A$.)
Then $T \circ S$ is a (compact, open) $G$-set containing
$(u,y)$.  For large $n$, $(u_n,y_n) \in T \circ S$.
Hence, for large $n$, there is $x_n \in X$ such that
$(u_n, x_n) \in T$ and $(x_n, y_n) \in S$.  The coordinate
projection maps are local homeomorphisms; consequently
$(u_n, x_n) \longrightarrow (u,x)$ and
$(x_n, y_n) \longrightarrow (x, y)$ in $P$. 
For all large $n$, $u_n$, $x_n$, and $y_n$ are
in the same equivalence class,
$u_n$ and $y_n$ are in $I$, and $x_n$
is in between $u_n$ and $y_n$;
so, $x_n \in I$.
 Thus, $(x_n, y_n) \in P \cap (I \bigtimes I)$
and hence $(x,y) \in \sc$, contradicting the assumption
that $(x,y) \in \s$.
\par
This proves that $(u,x) \in P, \; (x,y) \in \s  
\Longrightarrow (u,y) \in \s$.  The proof that
$(x,y) \in \s, \; (y,v) \in P \Longrightarrow
(x,v) \in \s$ is similar; the two implication together
show that $\s$ is an ideal set. 
\par
Next, we show that $\s$ is meet irreducible.  Suppose
that $\t_1$ and $\t_2$ are ideal sets and that
$\s = \t_1 \cap \t_2$.  Assume that $\s$ is a proper
subset of both $\t_1$ and $\t_2$.
\par
First observe that there is a point $(x,y) \in 
P \cap (I \bigtimes I)$ such that $(x,y) \in \t_1 \setminus
\s$.  Indeed, assume the contrary.  Since no point of  
$P \cap (I \bigtimes I)$ lies in $\s$,
 we have $P \cap (I \bigtimes I) 
\subseteq P \setminus \t_1$.  But $P \setminus \t_1$ is 
closed, so $\sc \subseteq P \setminus \t_1$.  This 
implies $\t_1 \subseteq \s$ (and therefore $\t_1 = \s$), 
contradicting our assumptions.
\par
Since $(x,y)$ is in $\t_1 \setminus \s$, we have
$(x,y) \in P \setminus \t_2$.  (Otherwise, $(x,y)
\in \t_1 \cap \t_2 = \s$, a contradiction.)
\par
Let $(a,b) \in P \cap (I \bigtimes I)$ and
let $u = \min \{a,x\}$ and $v = \max \{b,y\}$.
(Here min and max are with respect to the order
on $I$.)  Then we have
$$
(u,x) \in P, \; (x,y) \in \t_1, \; (y,v) \in P
\Longrightarrow (u,v) \in \t_1.
$$
Since $(u,v) \in P \cap (I \bigtimes I) \subseteq P \setminus \s$
we also have $(u,v) \notin \t_2$.  But since $(u,a) \in P$,
$(b,v) \in P$ and $\t_2$ is an ideal set, this implies that
$(a,b) \notin \t_2$.  As $(a,b)$ is arbitrary in 
$P \cap (I \bigtimes I)$, we obtain $P \cap (I \bigtimes I)
\subseteq P \setminus \t_2$.  Since the latter is a closed
set, this yields $\sc \subseteq P \setminus \t_2$, which implies
$\s = \t_2$, contrary to assumption.  This shows that
$\s$ must equal one of $\t_1$ or $\t_2$ and hence is
meet irreducible.
\qed \enddemo
\par
\remark{Remark} The mapping from intervals contained in some
equivalence class of $G$ to meet irreducible ideal sets
is not one-to-one, even in a context as simple as a 
refinement algebra.  Some meet irreducible ideal sets
can be written as the complement of $\sc$ for
a unique interval $I$ from a unique equivalence class.
For others, there is at least one such interval $I$
for each equivalence class from $G$.  It is also possible
that different intervals from the same equivalence class
yield the same meet irreducible ideal set.
(Here, the latitude lies in whether or not to include
``end points''.)
\endremark
\par
Theorem 2.1 has a converse, whose proof requires
the following elementary fact.

\proclaim{Fact 2.2} For each element $e$ of $T_n$, let
$\Id(e)$ denote the ideal in $T_n$ generated by $e$.
If $e_{ii}$, $e_{jj}$ and $e_{kk}$ are three diagonal
matrix units with $i < j < k$, then
$\Id(e_{ii}) \cap \Id(e_{kk}) \subseteq \Id(e_{jj})$.
\endproclaim

\proclaim{Theorem 2.3} 
With notation as above, let $\II$ be a meet irreducible 
ideal in $A$.
Then there is an interval $I$ contained in an equivalence 
class from $G$ so that the support set of $\II$ 
is $P \setminus \sc$.
\endproclaim 
\demo{Proof} The first step is to determine the equivalence
class which will contain $I$.  By Theorem 1.2, there is
a presentation
$$
\presen A A 1 2 3
$$
together with an MI-chain $(e_k)_{k \ge 1}$ for which
$\II$ is the largest ideal which contains no matrix 
unit $e_k$ from the MI-chain.  We shall use the MI-chain
to construct a decreasing sequence of projections
$p_1 \ge p_2 \ge p_3 \ge \dots $ with each
$p_k \in D_k = A_k \cap A_k^*$.  Each such decreasing
sequence of projections corresponds in a natural way
to a point of $X$ and thereby determines an equivalence
class in $G$.
\par
\comment
We need the following observation: for each
$i$, there is a matrix unit $f_{i+1}$ in $A_{i+1}$
such that $f_{i+1}^* e_{i+1}\phst e_{i+1}^* f_{i+1}\phst \le
e_i\phst e_i^*$.  To see this, recall that 
$e_{i+1} \in \Id_{i+1}(e_i)$.
Now $e_i$ is a sum of matrix units in $A_{i+1}$ and 
$\Id_{i+1}(e_i)$ is equal to the linear span of matrix units 
of the form $r e_{\tsub} s$, where $r$,\;$s$ and $e_{\tsub}$ are 
matrix units in $A_{i+1}$ and $e_{\tsub}$ is a subordinate  
of $e_i$.  
In particular, since $e_{i+1}$ is a matrix unit, 
there exist matrix units $f_{i+1}$ and $s$ in $A_{i+1}$ and 
a subordinate $e_{\tsub}$ of $e_i$ such that $e_{i+1} = f_{i+1} e_{\tsub} s$.
Since $e_{i+1} \ne 0$, we must have $f_{i+1}^* f_{i+1}\phst
= e_{\tsub}\phst e_{\tsub}^*$
 and $ss^* = e_{\tsub}^* e_{\tsub}\phst $.
We now have
$$
f_{i+1}^* e_{i+1}\phst e_{i+1}^* f_{i+1}\phst = 
 f_{i+1}^* f_{i+1}\phst e_{\tsub}\phst 
s s^* e_{\tsub}^* f_{i+1}^* f_{i+1}\phst
= e_{\tsub}\phst e_{\tsub}^* \le e_i\phst e_i^*,
$$
as desired.
\par
The sequence $p_n$ can now be defined:
$$
\aligned
p_1 &= e_1\phst e_1^*, \\
p_n &= f_2^* f_3^* \dots f_n^* e_n\phst e_n^* f_n\phst
 \dots f_3\phst f_2\phst,
\qquad n \ge 2.
\endaligned
$$
Let $z$ be the point in $X$ corresponding to the sequence
$(p_n)$.  The equivalence class of $z$, $\orb_z$, is
totally ordered by $P$.  
\endcomment
\par
Observe that $\Id_2(e_1)$ is equal to the linear span of
matrix units of the form $fsg$, where $f,s,g \in A_2$ and
$s$ is a subordinate of $e_1$ in $A_2$.  Since
$e_2$ is a matrix unit and is in $\Id_2(e_1)$, it has
this form.  In particular, there is a matrix unit $s_2$
in $A_2$ which is a subordinate of $e_1$ such that
$e_2 \in \Id_2(s_2)$.  Let $p_2$ and $q_2$ be the
range and domain projections of $s_2$;
 i.e., $p_2 = s_2^{\phst}s_2^*$ and
$q_2 = s_2^*s_2^{\phst}$.  If we let $p_1 = e_1^{\phst}e_1^*$
and $q_1 = e_1^*e_1^{\phst}$, then we have
$p_1 \ge p_2$ and $q_1 \ge q_2$.  Note also that
$e_2 \in \Id_2(p_2)$ and $e_2 \in \Id_2(q_2)$, since
both of these ideals contain $\Id_2(s_2)$.
\par
By property (B) for MI-chains, $e_3 \in \Id_3(e_2)$;
 consequently $e_3 \in \Id_3(s_2)$.
Therefore, there is a matrix unit $s_3$ in $A_3$
which is subordinate to $s_2$ (and hence to $e_1$)
such that $e_3 \in \Id_3(s_3)$.  Let
$p_3 = s_3^{\phst}s_3^*$ and $q_3 = s_3^*s_3^{\phst}$.
We have $p_2 \ge p_3$, $q_2 \ge q_3$, 
$e_3 \in \Id_3(p_3)$ and $e_3 \in \Id_3(q_3)$.
\par
It is now clear that an inductive argument will yield
a sequence of matrix units $s_n$ in $A_n$ with
range projections $p_n$ and domain projections $q_n$
such that:
 \item{1)} $s_1 = e_1$,
\item{2)} $s_{n+1}$ is a subordinate of $s_n$, for all $n$,
\item{3)} $e_n \in \Id_n(s_n)$, $e_n \in \Id_n(p_n)$
 and $e_n \in \Id_n(q_n)$, for all $n$, and
\item{4)} $p_n \ge p_{n+1}$ and $q_n \ge q_{n+1}$, for all $n$.
\fp
Clearly, $(p_n)$ and $(q_n)$ give points $p$ and $q$ in $X$.  
Since $(p,q) \in \ehat_1$, $p$ and $q$ determine the same 
equivalence class in $G$.  
This is the equivalence class which will contain $I$.
\par
If $x \in X$ then, for each $k$, there is a unique
minimal projection $x_k$ in $D_k$, the diagonal of $A_k$,
such that $x \in \xhat_k$. 
Define $I$ as follows:
$$
I = \{ x \in \orb_p \: e_k \in \Id_k(x_k)
  \text{ for all large } k\}.
$$
Note that both $p$ and $q$ are in $I$.
\par
For later use we need an observation.  Fix $k >1$.
Let $f_k$ and $g_k$ be matrix units in $A_k$
for which $e_k = f_k s_k g_k$.  For each $n \ge k$,
let $\tilde s_n = f_k s_n g_k$ and let $\tilde p_n$
and $\tilde q_n$ be the range and domain projections
of $\tilde s_n$.  
Then $\tilde s_n$, $\tilde p_n$, $\tilde q_n$, $n \ge k$ satisfy 
properties analogous to the properties 1)--4) above 
for $s_n$, $p_n$, $q_n$, $n \ge 1$.  
In particular, $e_n \in \Id_n(\tilde p_n)$ and
$e_n \in \Id_n(\tilde q_n)$ for all $n \ge k$;
the points $\tilde p$ and $\tilde q$ in $X$
corresponding to $(\tilde p_n)$ and $(\tilde q_n)$
lie in $I$; and $(\tilde p, \tilde q) \in \ehat_k$.
Thus, for any $e_k$ we can construct a point 
$(\tilde p, \tilde q)$ in $\ehat_k \cap (I \bigtimes I)$.
\par
We must show that $I$ is an interval in $\orb_z$.
Suppose $w < x < y$ where $w, y \in I$ and
$(w_k)$, $(y_k)$ are the nested sequences of projections
associated to $w$ and $y$.
Recall that we sometimes write $w \le x$ when $(w,x) \in P$.
There is an integer $N$ such that, for any $k \ge N$,
all of the following are true:
\item{i)} $e_k \in \Id_k(w_k)$,
\item{ii)} $e_k \in \Id_k(y_k)$,
\item{iii)} there is a matrix unit in $A_k$ with initial
projection $x_k$ and range projection $w_k$, and
\item{iv)} there is a matrix unit in $A_k$ with initial
projection $y_k$ and range projection $x_k$.
\fp
Now, $A_k$ is a maximal triangular subalgebra of a finite 
dimensional \cstar-algebra and so is a direct sum of
$T_n$'s.  Conditions iii) and iv) imply that $w_k$,
$x_k$ and $y_k$ all lie in the same summand; furthermore
within that summand $x_k$ lies in between $w_k$ and
$y_k$ in the diagonal ordering.  Since the context is
now that of a $T_n$, Fact 2.2 tells us that
$\Id_k(w_k) \cap \Id_k(y_k) \subseteq \Id_k(x_k)$.
In particular, $e_k \in \Id_k(x_k)$.  Since this
holds for any $k \ge N$, $x \in I$, this proves that $I$ 
is an interval.
\par
It remains to show that $\II$ has support set $P \setminus \sc$.  
Let $\II'$ be the ideal with support set $P \setminus \sc$.
\par
Suppose  $e$ is a matrix unit which is in $\II$ but
not in $\II'$.
Then $\ehat \cap (I \bigtimes I) \ne \varnothing$, so there are 
points $x,y \in I$ such that $(x,y) \in \ehat$.  
There is an integer $k$ such that $e \in A_k$, 
$e_k \in \Id_k(x_k)$, and $e_k \in \Id_k(y_k)$.  
If $f_k$ is the matrix unit in $A_k$ for which $(x,y) \in {\fhat}_k$, 
then $f_k$ is a subordinate of $e$.  
Since $f_k$ generates $\Id_k(x_k) \cap \Id_k(y_k)$, 
we have $e_k \in \Id_k(f_k)$.  
This implies that $f_k \notin \II$ and hence $e \notin \II$, 
a contradiction.  
Thus, $\II \subseteq \II'$.
\par
All that remains is to prove that $\II' \subseteq \II$.  
We observed earlier that, for each $k$, there is a 
point $(\tilde p, \tilde q) \in \ehat_k \cap (I \bigtimes I)$.  
Thus $e_k \notin \II'$, for all $k$.  
Suppose that $e$ is a matrix unit which is not in $\II$.
By the definition of $\II$, $e_k \in \Id_k(e)$ for some $k$.  
But this means that $e \notin \II$ for otherwise 
we would have $e_k \in \II'$, a contradiction.  
Thus $\II' \subseteq \II$. \qed
\enddemo
For each meet irreducible ideal, we can use the associated
interval $I$ to construct a nest representation 
whose kernel is the ideal.

\proclaim{Theorem 2.4} 
With notation as above, let $\II$ be a meet irreducible 
ideal in $A$
with associated interval $I$ as in Theorem 2.3.
Then there is a nest representation of $A$ acting on 
the Hilbert space $\hs$ whose kernel is $\II$.
\endproclaim
\demo{Proof}
Let $\{ \dx \: x \in I \}$ be the standard orthonormal basis 
for $\hs$.
Define $\pi$ on the matrix units in $A$ by, for a matrix 
unit $e$ and a basis vector $\dy$, setting
$$
\pi(e) \dy = \left\{
\alignedat2
&0, &&\qquad \text{if there is no } x \in I 
  \text{ such that } (x,y) \in \ehat, \\
&\dx, &&\qquad \text{if there is } x \in I 
  \text{ such that } (x,y) \in \ehat.
\endalignedat
\right.
$$
Since $\ehat$ is a $G$-set, $\pi(e) \dy$ is well defined.
Thus, $\pi(e)$ is a partial isometry in $\BB(\hs)$.
It is straightforward to check that 
$\pi(ef) \dy = \pi(e) \pi(f) \dy$ for any two
matrix units $e, f \in A$; 
so, the linear extension of $\pi$ to the algebra (not closed) generated 
by the matrix units of $A$ is an algebra homomorphism. 
\par
The obvious extension of $\pi$ to the matrix units of the
\cstar-envelope, $B$, of $A$ and the algebra generated
by these matrix units is also a $*$-algebra homomorphism.  
Furthermore, it has norm 1, since its restriction to each
$B_k = C^*(A_k)$ is a representation of a \cstar-algebra.
Since $\pi$ has norm 1, it extends to a representation
of $A$ acting on $\hs$.
\par
If $M$ is an invariant subspace for $\pi$ and if $\dy \in M$, 
then $\dx \in M$ for all $x \in I$ with $x \le y$.  
This is immediate, since $x \le y$ means that there is 
a matrix unit $e \in A$ with $(x,y) \in \ehat$.
Thus, if $M$ is an invariant subspace for $\pi$,
there is an initial segment $S$ of $I$ such that
$M$ is generated by $\{\dx \: x \in S \}$.  
This implies that the invariant subspaces for $\pi$ are 
totally ordered by inclusion.  
Thus, $\pi$ is a nest representation.
\par
Recall from Theorem 2.3 that $\II$ has support set $P \setminus \sc$.  
If $e$ is a matrix unit in $A$ then $\pi(e) = 0$ 
if, and only if, $\ehat \cap (I \bigtimes I) = \varnothing$.  
If $\ehat \cap (I \bigtimes I) = \varnothing$, then	
$P \cap (I \bigtimes I)$ is disjoint from the
open set $\ehat$; hence $\sc$ is disjoint from $\ehat$.  
Thus $\ehat \subseteq P \setminus \sc$ and so $e \in \II$.
Since ideals are generated by the matrix units which
they contain, it follows that $\ker \pi \subseteq \II$.
On the other hand, if $e$ is a matrix unit in $\II$,
then we have $\ehat \subseteq P \setminus \sc$, whence
$\ehat \cap (I \bigtimes I) = \varnothing$ and $e \in \ker\pi$.  
Thus $\II \subseteq \ker\pi$ and we have equality.\qed
\enddemo
\par
\head 3. Ideal Sets for Trivially Analytic Algebras
\endhead
\par
In this and the next section,
 we shall focus primarily on TAF algebras which are
analytic.  
An analytic subalgebra of an AF \cstar-algebra is
automatically strongly maximal triangular.  So the
results of the previous section apply in this	setting.
The class of trivially analytic subalgebras of AF \cstar-algebras
is fairly extensive; it includes, for example, all full nest
algebras.  These are algebras with a presentation of the
form
$$
\Tpresen {n_1} {n_2} {n_3} A  
$$
subject to the requirement that each embedding carries
the nest of invariant projections of $T_{n_i}$ into
the invariant projections of $T_{n_{i+1}}$.
The well-known refinement algebras form a subfamily of
the family of full nest algebras.
\par
In this section we shall give a complete description of
all the meet irreducible ideals in a trivially analytic
TAF algebra with an injective 0-cocycle
via a description of the meet irreducible	
ideal sets of the spectrum of the algebra.  This is the
setting most analogous to the finite dimensional context.
It is the context with the most intuitive picture
of meet irreducible ideal sets.
\par
\remark{Remark}
The description of the meet irreducible ideal sets is actually 
valid in a somewhat more general context, which we outline in
this remark.
The basic properties that we need for the description of 
the meet irreducible ideal sets are the following:
\item{1.} Each equivalence class from $G$ is countable
and dense in $X$.
\item{2.} The two projection maps from $X \bigtimes X$ 
to $X$ when restricted to $G$ are open and continuous
with respect to the groupoid topology on $G$.
\item{3.} There is a total order $\preceq$ on $X$ which,
 on each equivalence class from $G$,
agrees with the order 
induced by $P$.  Furthermore, the order topology on $X$ 
agrees with the original, Gelfand topology on $X$.
\fp The first of these properties implies that the 
groupoid \cstar-algebra associated with $G$ is simple.
The second property is equivalent to $G$ being $r$-discrete
and admitting a left Haar system. 
See \cite{R, Prop. 1.2.8}.
\par
The third property is the critical one for our purposes.
The existence of a total order on $X$ with these properties
follows immediately from the existence of a trivial
cocycle which is the coboundary of an injective 0-cocycle $b$:
define $x \preceq y$ iff $b(x) \le b(y)$.  
\par
The existence of a total order with property 3 is almost, but
not quite equivalent to the existence of a trivial cocycle
on $X$ which is the coboundary of an injective function.  
Equivalence requires one additional property: the
order $\preceq$ has at most countably many gaps.  (A gap is
a pair of elements from $X$
with no intermediate elements from $X$.)
\par	
It is not difficult to construct an example of a triple
$(X, P, G)$ which meets all of the properties above except
that it has uncountably many gaps with respect to the order
on $X$.  (Basically, construct a Cantor like set from the
interval $[0,1]$ doubling the irrational points instead
of the rational points. For the groupoid $G$ take the union
of all sets of the form $\{(qx,x) \: x \in A\}$, where $A$ is
some open interval from $X$ and $q$ is a positive rational
number with the property that $qX \subseteq X$.)  
The \cstar-algebra built on such a groupoid will be
inseparable and will fail to have most of the nice
properties that groupoid \cstar-algebras usually enjoy,
so this example is of dubious interest.)
\par
If $X$ does have countably many  gaps, construct a 
one-to-one, continuous map $b \: X \longrightarrow \Bbb R$
as follows.  Let $S$ be a countable dense subset of $X$ which
does not contain any points which have either an immediate
successor or an immediate predecessor.  Let 
$a \: S \longrightarrow [0,1]$
 be a monotonically increasing map of $S$
onto a countable, dense subset of $[0,1]$.  Extend $a$ to
a continuous map (also denoted by $a$) of $X$ onto $[0,1]$.
The map $a$ is increasing, but not
necessarily one-to-one.  In particular, if $x$ is an
immediate predecessor of $y$, then $a(x) = a(y)$.
Let $\{ (x_n, y_n) \}$ be an enumeration of all the gap
pairs from $X$.  For each $x$, let $\beta(x) =
\{n \: y_n \le x \}$.  Define $b \: X \longrightarrow
[0,2]$ by
$$
b(x) = a(x) + \sum_{n \in \beta(x)} \frac 1 {2^n} .
$$
The function $b$ has all the desired properties.
\par
The description of all the meet irreducible ideals
in a trivially analytic TAF algebra
with injective 0-cocycle can be verified
making use of only the properties of the spectral
triple listed above.  It is not necessary to use
the countability of the gap points nor the fact that
the enveloping \cstar-algebra is AF.  The argument, however, 
is long, tedious, and of little intrinsic interest.
Consequently, we will instead make use of Theorems~2.1 and~2.3
to provide a much more palatable verification at the 
expense of a slight loss of generality.  \qed
\endremark

For the following, assume that $A$ is a trivially 
analytic TAF algebra with diagonal $D$ and enveloping
\cstar-algebra $B$, which is simple.  
Let $(X, P, G)$ be the spectral
triple for $(D, A, B)$.  Let $\preceq$ be a total order
on $X$ which agrees with $P$ on equivalence classes from
$G$ and assume that the order topology agrees with the
original (Gelfand) topology on $X$.  If a point $a \in X$
has an immediate successor, we say that $a$ has a {\it gap
above}.  Similarly, if $b$ has an immediate predecessor,
then $b$ has a {\it gap below}.
\definition{Notation} For each pair of points $a, b \in X$
let
$$
\align
\sab &= \{(x,y) \in P \: x \prec a \text{ or } b \prec y \} \\
\tab &= \sab \cup \{(a,b)\}.
\endalign
$$
\enddefinition
\par
Observe that the set $\sab$ is an open subset of $P$ which
satisfies the ideal property.  Thus, it is always the support
set for an ideal in $A$.  The set $\tab$ also satisfies the
ideal property, but it need not be open.  
It will be an open subset of $P$ precisely when  $(a,b) \in P$ 
and there is a neighborhood, $\nu$, of $(a,b)$ 
such that $\nu \setminus \{(a,b)\} \subseteq \sab$.  When this
is the case, $\tab$  is an ideal set. 
In a refinement algebra, $\tab$ is an ideal set for all
$(a,b) \in P$.  In a full nest algebra, there may be
points $(a,b)$ for which $\tab $ is not open.
 In the following,
we shall always assume that $\tab$ is an ideal set.
\par
If $b \prec a$, then $\sab = P$.  If $a=b$, then
$\tab = P$ and $\sab$ is a maximal ideal (with codimension 1).
The ideal set $P$ is meet irreducible by default and each
$\s_{a,a}$ is trivially meet irreducible.  Consequently,
in the proof of Theorem 3.1, we always assume $a \prec b$.
\par
All meet irreducible ideal sets for a trivially analytic algebra
are described in the following theorem.
\proclaim{Theorem 3.1}  
Let $A$ be a trivially analytic TAF algebra whose spectral triple is
 $(X, P, G)$.  
Let $\preceq$ be a total order on $X$ compatible with the spectral triple.
The following is a complete list of all the meet irreducible
ideal sets in $P$:
\item{1.} $\sab$ if $(a,b)\in P$.
\item{2.} $\sab$ if $(a,b) \notin P$ and there is either 
no gap above for $a$ or no gap below for $b$.
\item{3.} $\tab$ if $(a,b) \in P$, there is either no gap 
above for $a$ or no gap below for $b$, and $\tab$ is open.
\endproclaim
\demo{Proof} Let $\s$ be a meet irreducible ideal
set contained in $P$.  By Theorem 2.3, there is an
equivalence class, $\orb_z$, from $G$ and
an interval $I \subseteq \orb_z$ such that
$\s = P \setminus \sc$.  Let $a = \inf I$ and
$b = \sup I$.  The inf and sup are taken in $X$
with respect to the order $\preceq$; the compactness
of $X$ guarantees that the inf and sup exist.
\par
We observe first that $\sab \subseteq \s$.
Indeed, suppose that $(x,y) \in P$ and $(x,y) \notin \s$.  
Then $(x,y) \in \sc$.  
Now, $\sc \subseteq P \cap (\overline I \bigtimes \overline I )$ 
(the containment may be proper), so $a \preceq x \preceq b$
and $a \preceq y \preceq b$.  But this shows that
$(x,y) \notin \sab$.  Thus, $\sab \subseteq \s$.

The next observation is that $\s \subseteq \tab$.
Indeed, suppose that $(x,y) \in P \setminus \tab$.  Then we
know that $a \preceq x,y \preceq b$ and $(x,y) \ne (a,b)$.
If both $x \ne a$ and $y \ne b$, then there is
an open neighborhood, $\nu$, of $(x,y)$ which is contained 
in $P \setminus \tab$.
We may further assume that the projection maps $\pi_1$ and
$\pi_2$ are homeomorphisms on $\nu$.  In particular,
$\pi_1(\nu)$ is an open set in $X$ which contains $a$.  
Consequently,
there is a point $u \in I \, \cap \pi_1(\nu)$.  It follows that
there is a point $(u,v) \in P \cap (I \bigtimes I)$ which
is in $\nu$.  This shows that $(x,y) \in \sc$.  If $x = a$ and
$a$ has no gap above, then $y \prec b$ and we may argue in
much the same way to conclude that $(a,y) \in \sc$.  If $x =a$
and $a$ has a gap above, then $a \in I$.  Since $y \prec b$,
we also have $y \in I$; in particular, $(a,y) \in \sc$.  The
case in which $y = b$ is handled in an analogous fashion. 
This proves that $\s \subseteq \tab$.

Since $\sab$ and $\tab$ differ by only one point, we have
shown that every meet irreducible ideal set has one of the
two forms $\sab$ or $\tab$.  To show that every meet
irreducible ideal set  is on the list in the theorem, 
we just have to show that the ideals of the form $\sab$
and $\tab$ which are not on the list are not meet
irreducible.

To that end, fix $(a,b) \in G \times G$ and let
$$
\align
\rho_1 &=\{ (x,y) \in P \: x \preceq a \text{ or } b \prec y\}  \\
\rho_2 &= \{ (x,y) \in P \: x \prec a \text{ or } b \preceq y \}.
\endalign
$$
Suppose that $(a,b) \notin P$ and that $a$
has a gap above and that $b$ has a gap below.
Since $a$ has a gap above and $b$ has a gap below, both
$\rho_1$ and $\rho_2$ are open and therefore ideal sets.
It is easy to check that $\sab$ is unequal to either 
$\rho_1$ or $\rho_2$ and that $\sab = \rho_1 \cap \rho_2$.
Thus $\sab$ is not meet irreducible when $(a,b) \notin P$,
$a$ has a gap above and $b$ has a gap below.

Suppose that $(a,b)\in P$ and $a$ has a gap above
and $b$ has a gap below.  
We  also assume that $a \ne b$, since otherwise $\tab = P$.  
Again, the hypotheses insure that $\rho_1$ and $\rho_2$ are ideal
sets which are unequal to $\tab$ and that $\tab = \rho_1 \cap \rho_2$.  
Thus, $\tab$ is not meet irreducible when $a$ has a gap above or 
$b$ has a gap below.

It remains only to show that the ideal sets on the list
are in fact meet irreducible.  This can be done by direct
argument or with the help of Theorem 2.1.  We will sketch
the argument which employs Theorem 2.1.
\par	
Suppose that $(a,b) \in P$.  Let $I =
\{ x \in \orb_a \: a \preceq x \preceq b \}$.
Then $\sab = P \setminus \sc$ and is therefore meet irreducible.
Note that this is the only choice for $I$ which works
in this case.  In subsequent cases the choice of $I$
will not be unique.
\par
Suppose that $(a, b) \notin P$ and $a$ does not have a gap
above.  In this case, let
$I = \{ x \in \orb_b \: a \prec x \preceq b \}$.
Then $\sc = \{ (x,y) \in P \: a \preceq x,y \preceq b \}$ and
$\sab = P \setminus \sc$ and so is meet irreducible.
\par
Suppose that $(a,b) \notin P$ and $b$ does not have a 
gap below.  This time we  
let $I = \{ x \in \orb_a \: a \preceq
x \prec b \}$.  Then $\sab = P \setminus \sc$ and is meet
irreducible.
\par
In the case in which  $a$ has no gap above 
and $b$ has no gap below, 
we take $I = \{ x \in \orb_z \: a \prec x \prec b \}$,
where $z$ is an arbitrary element of $X$.  Again
we get $\sab = P \setminus \sc$.
\par
Ideal sets of the form $\tab$ remain.
 Suppose that $(a,b) \in P$ and that $a$ has no gap above.
Let $I = \{ x \in \orb_b \: a \prec x \preceq b \}$.
Since  $\tab$ is an ideal set,  $(a,b)$
lies in an open neighborhood $N$ which is a subset of $\tab$
and therefore disjoint from $P \cap (I \bigtimes I)$.
This shows that $(a,b) \notin \sc$.  The rest of the argument
needed to show that $\tab = P \setminus \sc$ is similar to
what has been done before.  Thus $\tab$ is meet irreducible
when $a$ as no gap above.
\par
The argument that $\tab$ is meet irreducible when $b$ has
no gap below is analogous the the preceding one.  As in the
case for $\sab$, when neither $a$ has a gap above nor $b$ has
a gap below, there are many choices for the interval $I$ which
will yield $\tab = P \setminus \sc$.
\qed 
\enddemo
\head 4. Ideal Sets and the Extended Asymptotic Range
\endhead
In the first part of this section we gather some 
results about ideal sets
in the spectrum of a general analytic TAF algebra
whose enveloping \cstar-algebra $B$ is simple.
We then give some further results in the case in which
the extended asymptotic range of the cocycle (to be defined
below) is $\{0, \infty\}$.
Throughout this section $I$ will be an interval from
an equivalence class from $G$.  
As before, we write $x \le y$ when $(x,y) \in P$.  
We do not assume that there is an order on $X$ which extends $P$.  
The cocycle $c$ on $G$ will, in general, not be a coboundary.
The simplicity of $B$ is equivalent to the density
in $X$ of each equivalence class from $G$.
\par
\definition{Definitions} 
For any subset $E \subset X$, we say that $E$ is {\it increasing\/}
if $x \in E$ and $x \le y$ imply $y \in E$.  
We define {\it decreasing\/} in an analogous fashion.  
For an interval $I$, if the restriction of the cocycle $c$ to 
$I \bigtimes I$ is bounded, 
we say that $I$ is {\it finite with respect to\/} $c$.  
If $c|I \bigtimes I$ is unbounded, we say that $I$ is 
{\it infinite with respect to\/} $c$.
\enddefinition
\par 
Note that, unlike $I$, the set $E$ is not totally ordered 
by $\le$ (i.e. by $P$).
We also point out that infinite intervals exist only when 
the cocycle is not trivial.
(Trivial cocycles are necessarily bounded.)
\par
There is a considerable difference between the properties
of finite intervals and the properties of infinite intervals.
First we gather some results about infinite intervals.
We shall learn shortly that infinite intervals are of little
interest -- they yield only the trivial 0-ideal.
\par
\proclaim{Lemma 4.1} Suppose that $I$ is an interval from
an equivalence class which is infinite with respect to $c$.
Then $I$ is either increasing or decreasing.
\endproclaim
\demo{Proof} Let $\orb_a$ be the equivalence class
which contains $I$.  Assume that $I$ is neither increasing
nor decreasing.  Then there exist an element $y \in I$ 
and an element $z \in \orb_a$ such that $y < z$ and
$z \notin I$.  Also, there exist an element $x \in I$
and an element $w \in \orb_a$ such that $w < x$ and
$w \notin I$.  Since $I$ is an interval, no element
of $I$ can be less than $w$ nor greater than $z$.
Thus, if $(s,t) \in I \bigtimes I$, we have either
$w < s \le t < z$ or $w < t \le s < z$.  In particular,
the cocycle property implies that
$$
\align
0 \le c(s,t) &\le c(w,z) \qquad \text{if } s \le t, 
  \text{ and} \\
0 \le c(t,s) & \le c(w,z) \qquad \text{if } t \le s.
\endalign
$$
Thus, $|c(s,t)| \le c(w,z)$ in all cases and $I$ is
finite with respect to $c$ -- contrary to assumption.
This shows that $I$ is either increasing or decreasing. 
\qed
\enddemo
If $\nu$ is an open $G$-set contained in $G$, and if
$x \in \pi_1(\nu)$, then there is a unique element $y\in X$
such that $(x,y) \in \nu$;  
in this situation, we shall often write
$y  = \nu(x)$. 
We thereby identify $\nu$ with
a partial homeomorphism of $X$ into $X$.  In effect, we
are using the same symbol for the partial homeomorphism
and for its graph.  If $V$ is an open subset of $X$, we
let $\nu (V)$ denote $\{\nu(x) \: x \in V \cap \pi_1(\nu)\}$.
\par
\proclaim{Proposition 4.2} 
Let $I$ be an infinite interval from an equivalence class.  
Then $I$ is dense in $X$.
\endproclaim
\demo{Proof} Let $V = X \setminus \overline{I}$.  We have
to show that $V = \varnothing$.  Suppose that $V$ is not
empty.
We know from Lemma 4.1 that $I$ is either increasing or
decreasing.  Assume that it is increasing.  (If $I$ is
decreasing, a similar argument to the one below will also
yield a contradiction.)

 From the density of equivalence classes, it follows that
$X = \bigcup \nu(V)$, where the union is taken over all
compact, open $G$-sets $\nu$ (See [R]).  
However, $X$ is a compact set, so there are 
finitely many compact, open $G$-sets $\nu_1, \dots, \nu_k$ 
so that $X = \bigcup_{j=1}^k \nu_j(V)$.  
Since each $\nu_j$ is compact and $c$ is continuous, 
there is $M$ such that $c|\nu_j < M$, for all $j$.

Since $c$ is unbounded on $I$, there are points $t,x \in I$
such  that $c(t,x) >M$.  The $\nu_j(V)$ cover $X$, so
there is $j$ such that $x \in \nu_j(V)$; i.e., there is
$v \in V$ such that $x = \nu_j(v)$.  
Now $v \notin I$ (since $V$ is the complement of $\overline{I}$) 
and $t \in I$.  
The fact that $I$ is increasing implies that $v < t$.  
Also, since $(v,x)\in\nu_j$ and $c <M$ on $\nu_j$, we have $c(v,x) < M$.  
Thus we have $c(v,x) = c(v,t) + c(t,x)$ with $c(v,x) <M$ and $c(t,x) > M$.  
This implies that $c(v,t) < 0$.  
But that means that $t < v$, contradicting the observation above that $v < t$.
Thus we conclude that $V = \varnothing$ and $\overline{I} = X$. \qed
\enddemo

The following proposition is false without the assumption
that $I$ is either increasing or decreasing.

\proclaim{Proposition 4.3} 
Let $I$ be an interval contained in an equivalence class from $G$.  
Assume that $I$ is either increasing or decreasing and 
that $I$ is dense in $X$.  
Then $P \cap (I \bigtimes I)$ is dense in $P$.  
Consequently, the meet irreducible ideal associated with $I$ 
is the 0-ideal.
\endproclaim
\demo{Proof} We assume that $I$ is increasing.  The proof
when $I$ is decreasing is similar, as usual.
Let $(x,y) \in P$.  Let $\nu$ be an open $G$-set such
that $(x,y) \in \nu \subset P$ and the coordinate projections are
homeomorphisms on $\nu$.  Since $\overline{I} = X$, there is
a sequence $x_k \in \pi_1(\nu) \cap I$ such that $x_k \rightarrow 
x$ in $X$.  Let $y_k = \nu(x_k)$.  Since $(x_k, y_k) \in P$ 
and $I$ is increasing, $y_k \in I$ for all $k$.
The coordinate projections are homeomorphisms on $\nu$, so
$y_k \rightarrow y$ in $X$ and $(x_k, y_k) \rightarrow (x,y)$
in $P$.  Thus, $(x,y) \in \sc$ and $\sc = P$. \qed
\enddemo
\proclaim{Corollary 4.4} With the same assumptions as above,
$I \bigtimes I$ is dense in $G$.
\endproclaim
\demo{Proof} This follows from the fact that $G = P \cup P^{-1}$.
\qed
\enddemo
\proclaim{Corollary 4.5} If $I$ is an infinite interval with
respect to the cocycle $c$, then the meet irreducible ideal 
 associated with $I$ is $\{0\}$.
\endproclaim
\demo{Proof} Combine Lemma 4.1 and Propositions 4.2 and 4.3. \qed
\enddemo
\par
Assume that the cocycle $c$ is $\Bbb Z$-valued.  
The standard algebras provide a class of examples with
$\Bbb Z$-valued cocycles.  
See \cite{PW} for more on the relationship between 
$\Bbb Z$-valued cocycles and standard embeddings.
\par
Let $I$ be an interval which is finite with respect to
$c$.  We claim that in the $\Bbb Z$-valued cocycle case,
the interval $I$ is in fact a set with finite
cardinality.  Indeed, let $x$ be any element from $I$
and define a function $\phi \: I \longrightarrow \Bbb R$
by $\phi(y) = c(x,y)$.  Observe that $\phi$ is one-to-one.
(If $\phi(y_1) = \phi(y_2)$, the $c(x,y_1) = c(x,y_2)$
and hence $c(y_1,y_2) = c(y_1, x) + c(x, y_2) = -c(x, y_1)
+c(x, y_2) = 0 $.  Since $c^{-1}(\{0\})$ is the diagonal,
 $y_1 = y_2$.)
Thus, $\phi$ is a bounded, integer valued, one-to-one
map on $I$.  It is now immediate that $I$ is a finite set.
\par
If $I$ is a finite set, then of course $ \sc
 = P \cap (I \bigtimes I)$.  Thus the complement of
the ideal set $\s$ associated with $I$ is a finite
subset of $P$.  If $\II$ is the ideal corresponding
to $\s$ and $\pi$ is the nest representation given by
Theorem 2.3, then the construction of $\pi$ implies
that $\pi$ acts on a finite dimensional Hilbert space.  
Consequently, $\II$ has finite co-dimension in $A$.  
Thus, we have the following proposition:
\par
\proclaim{Proposition 4.6} Let $A$ be an analytic TAF
algebra whose \cstar-envelope is simple and which
has a $\Bbb Z$-valued cocycle.  Then any non-trivial
meet irreducible ideal in $A$ has finite co-dimension.
\endproclaim
\par
The results about the complement of the ideal set for
a meet irreducible ideal in an analytic algebra with
a $\Bbb Z$-valued cocycle can be extended in modified 
form to a broader class of algebras.  For this we
need the concept of {\it asymptotic range\/} from
\cite{R, Definition I.4.3} and a modification of
asymptotic range from \cite{S, p. 345}.
\par
\definition{Definitions}
If $c$ is a real valued cocycle, the {\it range\/}
of $c$ is $R(c) = \overline {c(G)}$.  If $U$ is a 
non-empty open subset of $X$, $c_U$ will denote
the restriction of $c$ to  $G \cap (U \bigtimes U)$.
The {\it asymptotic range\/} of $c$ is
$R_{\infty} (c) = \bigcap R(c_U)$, where the
union is taken over all non-empty open subsets
of $X$.  We say that $\infty$ is an {\it asymptotic
value\/} of $c$ if, for every $M>0$ and every non-empty 
open subset $U \subseteq X$, $R(c_U) \cap [M, \infty) \ne
\varnothing$.  Finally, we define the {\it extended 
asymptotic range\/} of $c$ to be
$$
\ear = \left\{
\alignedat 2
&R_{\infty}(c) &\qquad &\text{if } \infty \text{ is not
  an asymptotic value of } c, \\
&R_{\infty}(c) \cup \{\infty\} &\qquad & \text{if }
  \infty \text{ is an asymptotic value of } c.
\endalignedat  \right.
$$
\enddefinition
\par
It is shown in [S] that the extended asymptotic value is
an invariant for the algebra (with respect to isometric
isomorphism) and that there are only four
possible values for $\ear$: the sets $\{0\}$, 
$\{0, \infty\}$, $\Bbb R \cup \{\infty\}$ and
$\lambda \Bbb Z \cup \{\infty\}$ for some $\lambda \ne 0$.
The first case, $\ear = \{0\}$ occurs if, and only if,
the cocycle $c$ is trivial.  
On the other hand, the standard algebras satisfy
$\ear = \{0, \infty\}$.  
\proclaim{Theorem 4.7}
Assume that $A$ is an analytic TAF algebra whose cocycle $c$ 
has extended asymptotic range $\ear = \{0, \infty\}$.  
Assume also that the \cstar-envelope of $A$ is simple.  
Let $I$ be an interval from an equivalence class of $G$
which is finite with respect to $c$.  
Then $\overline{I}$ has empty interior.  
Consequently, $\sc$ has empty interior in $P$; i.e., 
the ideal set $\s$ corresponding to $I$ is dense in $P$.
\endproclaim
\demo{Proof} Suppose $I \subseteq \orb_a$ and
that $\overline{I}$ has
non-empty interior.  We first observe that we
may as well assume, without loss of generality,
that $\overline{I} = X$.  Indeed, if 
the interior of $\overline{I}$ is non-empty,
 then there is
a compact open subset $V \subseteq X$ such that
$V \subseteq \overline{I}$.
  We can then simply replace $G$ by
$G$ restricted to $V$.  We just need to note that
$\tilde{R}_{\infty}(c|G \cap (V \bigtimes V))
= \ear = \{0,\infty\}$.
\par
The assumption that $I$ is finite with respect to
$c$ implies that there is a number $M$ such that
$|c(x,y)| \le M$ for all $x,y \in I$.
Since $\ear = \{0, \infty \}$, for every $x \in X$ and
every $\epsilon > 0$ we can find an open set
$U(\epsilon, x)$ containing $x$ such that
$$
R\bigl(c\,|\, G \cap (U(\epsilon,x) \bigtimes U(\epsilon,x))\bigr)
 \cap (\epsilon, 2M) = \varnothing.
$$
Consequently
$$
R\bigl(c\,|\,I \cap (U(\epsilon,x) \bigtimes U(\epsilon, x))\bigr)
\subseteq [- \epsilon, \epsilon].
$$
\par
Suppose that $x_n \in I$ and $x_n \rightarrow x$.
Then there is $N$ such that for $n \ge N$,
$x_n \in U(\epsilon, x)$.  Hence, for $n,m \ge N$,
$$
|c(a, x_n) - c(a, x_m)| = |c(x_m, x_n)| \le \epsilon.
$$
It follows that $(c(a,x_n))$ is a Cauchy sequence
and therefore has a limit.
\par
If $y_n$ is another sequence from $I$ such that
$y_n \rightarrow x$, then by the same argument
the ``interwoven'' sequence $c(a, x_1)$, $c(a, y_1)$,
$c(a, x_2)$, $c(a, y_2),\; \dots$ is also Cauchy.
Consequently, $\ds \lim_{n \to \infty} c(a, x_n)
= \lim_{n \to \infty} c(a, y_n)$.
\par
We now define
$$
g(x) = \lim_{n \to \infty} c(a, x_n), \quad 
\text{where } x_n \in I \text{ and } x_n \rightarrow x.
$$
The argument above shows that $g$ is well defined; we next 
show that $g$ is continuous.  Fix $x \in X$ and $\epsilon >0$.
We shall show that for any $y \in U(\epsilon, x)$, 
$|g(y) - g(x)| \leq \epsilon$, thereby verifying that $g$ is
continuous.
There exist sequences $x_n \in U(\epsilon, x)$
and $y_n \in U(\epsilon, x)$ such that $x_n \rightarrow x$ and
$y_n \rightarrow y$.  
But then $|c(a, x_n) - c(a, y_n)| = |c(y_n, x_n)| \le \epsilon$.  
This holds for all $n$, so $|g(x) - g(y)| \le \epsilon$.
\par
Since $g(x) = c(a, x)$ for all $x \in I$, it follows that
$c(x,y) = g(y) - g(x)$ on $I \bigtimes I$.
\par
By assumption, the cocycle $c$ is unbounded on $G$.
Since $G \cap (\orb_a \bigtimes \orb_a)$ is dense in $G$,
it follows that $c$ is unbounded on $\orb_a \bigtimes \orb_a$.
\par
Write $ \orb_a = \bigcup_{n=1}^{\infty} I_n$, where
 $I \subset I_1 \subset I_2 \subset \dots $ are intervals
in $\orb_a$ and the restriction of $c$ to $I_n \bigtimes I_n$
is bounded for every $n$.  By the arguments above, for each
$n$ there is a continuous function $g_n$ defined on $X$ such
that $c(x,y) = g_n(y) - g_n(x)$ on $I_n \bigtimes I_n$.
 From the definition of the $g_n$, it follows that, for
$m < n$, $g_n | I_m = g_m | I_m$.  This shows that
$g_n = g_m$ for all $n,m$ (since $g_n$ and $g_m$ are
continuous and $\overline{I_n} = \overline{I_m} =X)$.
\par
In particular, for $x \in I_n$, $g(x) = g_n(x) = 
c(a,x)$.  This holds for all $n$, so in fact
$g(x) = c(a,x)$ for all $x \in \orb_a$.
Therefore $c(x,y) = g(y) - g(x)$ on $\orb_a \bigtimes
\orb_a$.  By the continuity of $c$ and $g$ and the
density of $\orb_a \bigtimes \orb_a$ in $G$, we have
$c(x,y) = g(y) - g(x)$ for all $(x,y) \in G$.  But this
says that $c$ is a trivial cocycle and hence that
$\ear = \{0\}$, contrary to assumption. 
\par
We have now proved that $\overline{I}$ has empty
interior in $X$.  It follows that $P \cap 
\overline{I} \bigtimes \overline{I}$ has empty
interior in $P$.  Since $\sc \subset
P \cap \overline{I} \bigtimes \overline{I}$, 
the remaining assertions of the theorem follow.
\qed
\enddemo
\remark{Remark} If $\ear = \{0, \infty \}$, then
every interval contained in an equivalence class
is either dense in $X$ or else nowhere dense.  
It follows that if $\II$ is a proper meet irreducible
ideal in $A$ and if $e$ is any matrix unit from
$A$, then there is a diagonal projection $q$ such
that $qe \in \II$.
\endremark
\comment
\head Section ?
\endhead
\par
\proclaim{Lemma ?.1} Let $I$ be an interval
 contained in an equivalence class from $G$
and let $\II$ be the associated meet irreducible ideal.
 Let $\s = P \setminus \sc$.  Also let $s = 
 \sup \{c(x,y) \: (x,y) \in I \bigtimes I \}$ and
$t = \inf \{c(x,y) \: (x,y) \in \s \}$.
If $s = \infty$, then $\II = (0)$.  Furthermore
$c^{-1}(s, \infty) \subseteq \s \subseteq c^{-1} [t, \infty)$.
\endproclaim
\demo{Proof} The first assertion is a rephrasing of
corollary 4.5.  The second is immediate from the
observation $s = \sup c\,|\, I \bigtimes I =
\sup c\,|\, P \cap (I \bigtimes I) = \sup c \,|\, \sc
= \sup c\,|\,(P \setminus \s)$. 
\enddemo
\proclaim{Proposition ?.2} Add to the hypotheses of lemma
?.1 the assumption that the extended asymptotic range
 $\ear = \Bbb R \cup \{\infty\}$.
Then the interior of $c^{-1}(t, \infty)
\setminus \s$ is $\varnothing$.
\endproclaim
\demo{Proof} Assume that $c^{-1}(t, \infty)
\setminus \s$ has non-empty
interior.  Then there is a an open $G$-set $\nu \subseteq
c^{-1}(t, \infty) \setminus \s$.  By making $\nu$ smaller,
if necessary, we may assume that $\nu \subseteq
c^{-1}(t_1, \infty) \setminus \s$ where $t_1 > t$.
Fix $t_2$ such that $t < t_2 < t_1$.  Since $t < t_2$,
we can find an open $G$-set $\nu_1 \subseteq \s$ such
that on $\nu_1$ we have $t \le c < t_2$.
Since $\ear = \Bbb R \cup \{\infty\}$, 
$c^{-1}([0, t_1 -t_2) \cup (\pi_1(\nu) \bigtimes \pi_1(\nu_1)))
\ne \varnothing$. 
\par
 Let $(w_1, w_2) \in
c^{-1}([0, t_2 - t_1) \cap (\pi_1(\nu) \bigtimes \pi_1(\nu)))$.
Then
$$
c(w_1, \nu_1(w_2)) = c(w_1, w_2) + c(w_2, \nu_1(w_2))
< t_1 - t_2 +t_2 = t_1
$$
  while $c(w_1, \nu(w_1) \ge t_1$.
Hence $(\nu_1(w_2), \nu(w_1)) \in P$ and
$(w_1, \nu(w_1)) = (w_1, w_2) \circ (w_2, \nu_1(w_2)) \circ
(\nu_1(w_2), \nu(w_1)) \in P \circ \s \circ P \subseteq \s$.
  But this
contradicts the fact that $\nu \cap \s = \varnothing$. \qed
\enddemo
\endcomment

\head 5. Completely Meet Irreducible Ideals \endhead

The paper \cite{DH} studies strongly maximal TAF algebras with isomorphic 
lattices of ideals and the extent to which one can conclude that the
algebras (or appropriate subalgebras thereof) are isomorphic
or anti-isomorphic.  
The essential tool for this study, $\MIC(A)$, turns out to be 
equivalent to the family of all completely meet irreducible ideals of $A$.  
In this section we give a theorem for completely meet irreducible ideals 
analogous to Theorem 1.2 and use this theorem to establish the connection
with $\MIC(A)$.   
In the case in which an algebra is generated by its order preserving 
normalizer, there is a natural bijection between the spectrum of the 
algebra and the family of completely meet irreducible ideals.
\definition{Definition} An ideal $\II$ is said to be {\it completely
meet irreducible\/} provided that, 
whenever $\II = \bigcap_{\lambda \in \Lambda}
\II_{\lambda}$, we have $\II = \II_{\mu}$, for some $\mu \in \Lambda$.
\enddefinition
\definition{Definition} A sequence $(e_k)_{k \ge N}$ of matrix
units from $A$ will be called a {\it CMI-chain\/} if the following
three conditions are satisfied for all $k \ge N$:
\item{(A)} $e_k \in A_k$.
\item{(B)} $e_{k+1} \in  \Id_{k+1}(e_k)$.
\item{(C)} The ideal in $A$ generated by $e_k - e_{k+1}$ does 
not contain $e_j$, for any $j \ge N$.
\enddefinition
\par
Of course, conditions (A) and (B) are just the conditions for
the sequence to be an MI-chain.
\remark{Remark} The three conditions above imply that $e_{k+1}$
is a subordinate of $e_k$, for each $k$.  Indeed, suppose that
$e_{k+1}$ is not a subordinate of $e_k$.  Let $f$ be any subordinate
of $e_k$ and let $r = ff^*$ be the range projection of $f$ and
$s = f^*f$ the initial projection of $f$.  Observe that
$f = re_ks \in \Id_{k+1}(e_k)$ and $re_{k+1}s =0$.  Since
$e_{k+1}$ is not in the ideal in $A$ generated by $e_k - e_{k+1}$
and $f = r(e_k - e_{k+1})s$, we conclude that $e_{k+1}$ is not
in the ideal in $A$ generated by $f$.  In particular,
$e_{k+1} \notin \Id_{k+1}(f)$.  Since this is true for each
subordinate of $e_k$, it follows that $e_{k+1} \notin \Id_{k+1}(e_k)$.
But this contradicts condition $(B)$.
\endremark
\proclaim{Theorem 5.1} 
Let $A$ be a strongly maximal TAF algebra with some presentation.
If $\II$ is an ideal in $A$, then $\II$ is completely meet irreducible 
if, and only if, $\II$ is the ideal corresponding to a 
CMI-chain of matrix units in the presentation.
\endproclaim
\demo{Proof} Suppose that $\II$ is completely meet irreducible.
If $\JJ$ is the intersection of all ideals of $A$ that properly
contain $\II$, then, by hypothesis, $\JJ$ properly contains $\II$.
By the choice of $\JJ$, there is no ideal, $\KK$, such that 
$\II \varsubsetneq \KK \varsubsetneq \JJ$. 
\par
If $\presen A A 1 2 3$ is a presentation for $A$, then there is
some $N$ such that $\II \cap A_N \ne \JJ \cap A_N$.  It is
easy to see that for each $j \ge N$, there is exactly one
matrix unit, call it $e_j$, that is in $\JJ \cap A_j$ but
not in $\II \cap A_j$.  By construction, $\II$ is the largest
ideal that does not contain $e_j$ for any $j \ge N$.  Also,
$e_j - e_{j+1}$ must be in $\II$, since it is in $\JJ \cap A_{j+1}$
and $e_{j+1}$ is not a subordinate of it.  
Thus, $(e_j)_{j \ge N}$ is a CMI-chain.
(Note: this can also be proved by appealing to \cite{DH, Lemma 2}.)
\par
Suppose that $\II$ is the ideal corresponding to a CMI-chain
$(e_k)_{k \ge N}$.  First observe that condition (C) implies
that $e_k - e_{k+1} \in \II$, for each $k$.  Let $\JJ$ be an
ideal that properly contains $\II$.  By the definition of
$\II$, there is some $j \ge N$ so that $e_{j+1} \in \JJ$.
But since $e_j - e_{j+1} \in \II \subseteq \JJ$, this implies
that $e_j \in \JJ$.  It follows that $\JJ$ contains $e_N$.
Hence, for a set of ideals each of which properly contains $\II$,
the intersection will contain $e_N$, and thus will properly
contain $\II$.  This shows that $\II$ is completely meet
irreducible. \qed
\enddemo
\par
In order to show the connection between completely meet irreducible 
ideals and the theory in \cite{DH}, we need some definitions and
notation from Section 2 of that paper.
\definition{Definition} If $\II$ and $\JJ$ are ideals in $A$, we
call $[\II, \JJ]$ a {\it minimal interval\/} if $\II \subsetneq
\JJ$ and if, whenever $\KK$ is an ideal in $A$ with
$\II \subseteq \KK \subseteq \JJ$, then either $\KK = \II$ 
or $\KK = \JJ$.
\enddefinition
\definition{Definition} If $[\II, \JJ]$ is a minimal interval,
its cone is the set $\{\KK \: \JJ \subseteq \KK \vee \II \}$.
Let $\MIC(A)$ denote the set of all equivalence classes of
minimal intervals under the equivalence relation of equal cones.
\enddefinition
\remark{Remark} Each equivalence class of minimal intervals
contains a maximal representative, $[\II, \JJ]$.  Just take
$\II = \bigvee \II_{\lambda}$ and $\JJ = \bigvee \JJ_{\lambda}$,
where both spans are taken over all minimal intervals
$[\II_{\lambda}, \JJ_{\lambda}]$ in the equivalence class.
Thus, we could equally well define $\MIC(A)$ to be the set
of all maximal representatives. 
\endremark
\proclaim{Proposition 5.2} 
Let $A$ be a strongly maximal TAF algebra.  
The set of all completely meet irreducible ideals in $A$ 
coincides with $\MIC(A)$.
\endproclaim
\demo{Proof} Let $\II$ be completely meet irreducible.  Set
$\JJ$ equal to the intersection of all ideals in $A$ which
properly contain $\II$.  By the complete meet irreducibility
of $\II$, $\JJ$ properly contains $\II$.  So, $[\II, \JJ]$ is
a minimal interval and hence gives an element of $\MIC(A)$.
\par
Fix an element of $\MIC(A)$ and let $[\II, \JJ]$ 
be the maximal representative of the equivalence class.
If $[\II', \JJ']$ is any element of the equivalence
class, then $\II$ is the join of all ideals $\KK$
for which $\KK \vee \II'$ does not contain $\JJ'$
and $\JJ$ is the join of $\JJ'$ and $\II$.
Repeat the argument of Proposition 5.1 or invoke
\cite{DH, Lemma 2} to see that there is a CMI-chain
of matrix units $(e_j)_{j \ge N}$ for which $\II$ is
the associated ideal.  By Proposition 5.1, $\II$ is
completely meet irreducible. 
\qed \enddemo
\par
The final result of this section establishes a bijection
between the spectrum of an algebra and the set of
completely meet irreducible ideals, provided that the
algebra is generated by its order preserving normalizer.
Again, we need a few definitions.  
The first is from \cite{PPW}.
\par
\definition{Definition} The {\it diagonal order\/} is the
partial order defined on the collection of all projections
in the diagonal, $D$, of $A$ as follows:  $e \preceq f$ if
there is a normalizing partial isometry, $w$, in $A$ such
that $e = ww^*$ and $f = w^*w$.
\enddefinition
\definition{Definition} 
If $w$ is a normalizing partial isometry, then the map, 
$x \longrightarrow w^*xw$, induces a bijection between 
the diagonal projections which are subprojections of $ww^*$ 
and the diagonal projections which are subprojections of $w^*w$.  
We say that $w$ is {\it order preserving\/} if this map preserves 
the diagonal order restricted to its range and domain.  
We define the {\it order preserving normalizer} of $A$ to be 
the set of all normalizing partial isometries which are order
preserving.
\enddefinition
\remark{Remark} 
If $\tau$ is the graph of an order preserving partial isometry 
(i.e., $\tau$ is an order preserving $G$-set), then there cannot be 
distinct points $(x,y) \in \tau$ and $(u,v) \in \tau$ such that
$x \le u \le v \le y$, where, as usual, $x \le u$ means 
$(x,u) \in P$.
This can be easily seen  by looking at the action of $\tau$ on the 
sequences of diagonal matrix units which correspond to the points 
$x$,~$u$,~$v$,~$y$ in $X$.  
As in the previous section, when $(x,y) \in \tau$ and $(u,v) \in \tau$, 
we write $x = \tau(y)$ and $u = \tau(v)$;
thus $\tau$ order preserving says that we cannot have
$v \le y$ and $\tau(y) \le \tau(v)$.
\endremark
\par
The concept of an order preserving normalizer first appeared
in \cite{MS1} in a groupoid context; the term used there for
the graph of an order preserving normalizing partial isometry is
{\it monotone G-set\/}.  The order preserving normalizer
was studied by Power in \cite{P3}, where it was called
the {\it strong normalizer\/}.  Note that a sum of 
order preserving elements which is again a partial isometry
is order preserving if, and only if, the ideal generated by 
each summand contains none of the other summands.
\par
For the remainder of this section, we assume that the 
algebra $A$ is generated by its order preserving
normalizer.  Algebras with this property were characterized
in terms of their presentations in \cite{DHo}.  The characterization
involves embeddings which are locally order preserving.
\definition{Definition} Let $A_1$ and $A_2$ be triangular 
subalgebras of finite dimensional \cstar-algebras.  An
embedding $\phi\:A_1 \longrightarrow A_2$ is {\it locally
order preserving\/} if $\phi(e)$ is order preserving for
each matrix unit $e \in A$.
\enddefinition
An algebra $A$ is generated by its order preserving
normalizer if, and only if, there is a presentation
for $A$ such that for any contraction of the presentation,
the embeddings in the contraction are locally order 
preserving \cite{DHo, Theorem 18}.  Another way to put
this is that each matrix unit in $A_j$ is order
preserving in $A_k$ when it is viewed as an element
of $A_k$, for any $k >j$.  This is, of course, equivalent
to saying that there is a system of matrix units
such that each matrix unit is an order preserving 
partial isometry in $A$.
\proclaim{Theorem 5.3} 
Let $A$ be a strongly maximal TAF algebra which is 
generated by its order preserving normalizer.  
Then there is a bijection between the spectrum, $P$, of $A$ 
and the set of completely meet irreducible ideals in $A$.
\endproclaim

Although one way to prove this theorem is to combine Proposition 5.2 
and \cite{DH, Theorem 7}, we give two self-contained proofs.
All three arguments use essentially the same underlying map, 
but the first proof below uses the inductive limit structure 
while the second uses the groupoid structure.
In particular, the second proof is not limited to subalgebras
of AF \cstar-algebras.

\demo{Proof 1} 
Fix a presentation for $A$ with the property that every embedding 
is locally order preserving.  
By Proposition 5.1, there is a one-to-one correspondence between
completely meet irreducible ideals in $A$ and CMI-chains.
But when every embedding is locally order preserving, the
CMI conditions are satisfied by every chain $(e_j)$ for
which each $e_{j+1}$ is a subordinate of $e_j$.  
The proof is completed by observing that the collection of all
such chains is in natural one-to-one correspondence with
the spectrum, $P$, of $A$. \qed
\enddemo
\demo{Proof 2} 
Given $(x,y) \in P$, let $I = [x,y] =\{u \: x \le u \le y \}$ be a 
closed interval in an equivalence class and let $Q(x,y) = \sc$.  
Let $J(x,y)$ be the (meet irreducible) ideal whose support
is $P \setminus \sc$.  
We shall show that the map $(x,y) \longrightarrow J(x,y)$ is a 
bijection from $P$ onto the collection of completely meet irreducible 
ideals.
\par 
First, we make a useful observation.  
If $\tau$ is an order preserving $G$-set which contains $(x,y)$ 
then $\tau \cap Q(x,y) = \{(x,y)\}$.  
To see this, first note that 
$\tau \cap P \cap (I \bigtimes I) = \{(x,y)\}$ --- this is
just the remark after the definition of order preserving
partial isometry.  Secondly, if $(w,z) \in \tau \cap Q(x,y)$,
then there is a sequence $(x_n, y_n) \in P \cap (I \bigtimes
 I) $ such that $(x_n, y_n) \longrightarrow (u,v)$.  For large
$n$, $(x_n, y_n) \in \tau$;
 therefore, $(x_n, y_n) \in \tau \cap Q(x,y)$.
Thus, $x_n = x$ and $y_n = y$ for large $n$; this shows that
$u=x$ and $v=y$, verifying the observation.
\par
Next we show that each ideal, $J(x,y)$, is completely
meet irreducible.  
It is convenient to work with the complements of ideal sets, 
so suppose that $Q(x,y) = \overline{\bigcup F_{\alpha}}$, 
where each $F_{\alpha}$ is the complement in $P$ of an ideal set.  
Since $I = [x,y]$ is a closed interval, 
$(x,y) \in Q(x,y) = \overline{\bigcup F_{\alpha}}$; 
hence, there is a sequence of points $(x_n, y_n) \in F_{\alpha_n}$
such that $(x_n, y_n) \longrightarrow (x,y)$.  
Let $\tau$ be an order preserving $G$-set which contains $(x,y)$.
Since $\tau$ is open, there exists $k$ such that 
$(x_k, y_k) \in \tau$.  
Thus, $(x_k, y_k) \in \tau \cap Q(x,y) = \{(x,y)\}$.  
So $x_k = x$ and $y_k = y$ and, hence, $F_{\alpha_k} = Q(x,y)$.
\par
To see that the map $(x,y) \longrightarrow J(x,y)$ is
onto the family of completely meet irreducible ideals,
let $\II$ be such an ideal.  Let $\s$ be the support
set for $\II$.  Observe that if $(x,y) \in P \setminus \s$, 
then $Q(x,y) \subseteq P \setminus \s$.  Thus
$$ \bigcup_{(x,y) \notin \s} Q(x,y) = P \setminus \s.$$
(Equality follows from the fact that each $(x,y) \in Q(x,y)$.)
Thus $ \s = \bigcap_{(x,y) \notin \s} P \setminus Q(x,y)$
and hence $\II = \bigcap_{(x,y) \notin \s} J(x,y)$.  
By the complete irreducibility of $\II$, $\II = J(x,y)$, 
for some $(x,y)$.
\par
It remains to show that the mapping is one-to-one.  Assume
that $J(x,y) = J(u,v)$ for points $(x,y), (u,v) \in P$.
Then $Q(x,y) = Q(u,v)$.  Let $\tau_1$ and $\tau_2$ be
order preserving $G$-sets such that $(x,y) \in \tau_1$
and $(u,v) \in \tau_2$.  Since $(x,y) \in Q(u,v)$,
there is a sequence $(x_n, y_n) \in \tau_1$ such that
$(x_n, y_n) \longrightarrow (x,y)$ and 
$u \le x_n \le y_n \le v$.
For every $n$, $(x_n, y_n) \in Q(u,v)$; hence
$$
\bigcup_n Q(x_n, y_n) \subseteq Q(u,v) = Q(x,y).
$$
Since $ (x,y) = \lim_n (x_n, y_n) \in
\overline{\bigcup Q(x_n, y_n)}$, we have 
$Q(x,y) \subseteq \overline{\bigcup Q(x_n, y_n)}$.
Thus $Q(x,y) = Q(u,v) = \overline{\bigcup Q(x_n, y_n)}$.
Since $Q(x,y)$ is completely meet irreducible, there is
 $m$ such that $Q(x,y) = Q(u,v) = Q(x_m, y_m)$.
We have $u \le x_m \le y_m \le v$ and, also,
$(u,v) \in Q(x_m, y_m)$; hence, there are $z_k$, $w_k$
such that $(z_k, w_k) \longrightarrow (u,v)$ 
and $x_m \le z_k \le w_k \le y_m$, for all $k$.  Without
loss of generality, we may assume that $(z_k, w_k) \in 
\tau_2$, for all $k$.   But then 
$u \le x_m \le z_k \le w_k \le y_m \le v$. 
Since $\tau_2$ is order
preserving, we must have $u = z_k$, for every $k$.
Thus, $u \le x_m \le u$; i.e., $u = x_m$.  
We can replace the sequence $\{(x_n, y_n)\}_{n=1}^{\infty}$
by $\{(x_n, y_n)\}_{n=N}^{\infty}$ for every $N \in \Bbb N$;
hence we can find a subsequence $(x_{m_k})$ such that
$x_{m_k} = u$, for all $k$.  Therefore, $x = \lim x_{m_k}
=u$.    This shows that $(x,y) = (u,v)$ and the mapping
is one-to-one.\qed 
\enddemo

\head 6. A Distance Formula \endhead
\par
In this section we prove a distance formula for ideals in 
strongly maximal TAF algebras which is analogous 
to the distance formula for a nest algebra.
First we prove the distance formula for the special case of an 
elementary groupoid of type $n$ \cite{R, III.1.1}, 
i.e., the groupoid corresponding to $M_n(C(X))$ 
where $X$ is a suitable topological space.
Recall that we use $[i,j]$ for the set $\{ i,i+1,\ldots, j\}$.

\proclaim{Proposition 6.1} 
Let $X$ be a locally compact, second countable Hausdorff topological space, 
let $H = X \bigtimes [1,n] \bigtimes [1,n]$ and 
suppose that $Y \subseteq H$ satisfies $(x, (i,j)) \in Y$ implies 
$(x, (i',j')) \in Y$ for all $i'$, $j'$ with $i' \le i$, $j' \ge j$.
If $f \in C(H)$ satisfies
$$ \sup \left\{ \left\| 
    f|_{\displaystyle \{x\} \bigtimes [i_0, n] \bigtimes [1, j_0]} \right\|
    \: \{x\} \bigtimes [i_0, n] \bigtimes [1, j_0] \subseteq Y \right\}\le 1,
    \tag {$*$}
$$
where the norm is the matrix norm of the restriction of $f$, then
there is $g \in C(H)$ so that $g=f$ on $Y$ and, for each $x \in X$,
$$ 
   \left\| g|_{\displaystyle\{x\}\bigtimes [1,n]\bigtimes [1,n]} \right\|\le 1.
$$
\endproclaim
\demo{Proof} 
First we order the $n^2$ coordinates of $[1,n] \bigtimes [1,n]$
in such a way that $(n,1)$ is first, $(1,n)$ is last and, if 
$i_1 \ge i_2$ and $j_1 \le j_2$, then $(i_1, j_1)$ precedes $(i_2, j_2)$.  
There are clearly many ways to do this.
Write $Z_m$ for the first $m$ coordinates in this ordering.
Let $g_0 = f$.  
We define, inductively, $g_m \in C(H)$ so that
\item{(1)} $g_m = f$ on $Y$, and
\item{(2)} condition ($*$) is satisfied for $g_m$
in place of $f$ and $Y \cup (X \bigtimes Z_m)$ in place of $Y$.
\fp
Setting $g = g_{m^2}$ then completes the proof.
\par 
We start by defining, for $a \ge 0$ and $b \in \Bbb C$,
$$
h(a,b) = \left\{
\aligned
& 0, \\
& \frac b {|b|} \min (|b|, a), 
\endaligned
\right.
\quad
\aligned
& \text{if } b = 0, \\
& \vphantom{\frac b {|b|}} \text{if } b \ne 0.
\endaligned
$$
We have the following three properties: 
(a) $|h(a,b)| \le a$,
(b) if $|b| \le a$, then $h(a,b) = b$, and
(c) for continuous functions $a(x)$, $b(x)$ with $a(x)\ge 0$, $b(x)\in\Bbb C$, 
    the map $x \mapsto h((a(x), b(x))$ is continuous.
\par
By (c), we have $g_1 \in C(H)$ where $g_1$ is defined by
$$
g_1(x, (i,j)) = \left\{
\aligned
&h(1, f(x,(n,1))), \\
&f(x, (i,j)),
\endaligned
\right. \quad
\aligned
&\text{if } (i,j) = (n,1), \\
&\text{if } (i,j) \ne (n,1).
\endaligned
$$
Also, if $(x, (n,1)) \in Y$, then ($*$) above implies that
$|f(x,(n,1))| \le 1$ and (b) shows that $g_1(x,(n,1)) = f(x,(n,1))$.  
Hence we get $g_1 = f$ on $Y$.
If $(x,(n,1)) \notin Y$, then we get $|g_1(x, (n,1))| \le 1$
(by (a)) and, thus, ($*$) holds for $g_1$ and 
$Y \cup \{(x,(n,1)) \: x \in X\}$.  
This completes the initial induction step.
\par
Assume that $g_1, \dots, g_{m-1}$ are defined satisfying
(1) and (2).  To define $g_m$ we change $g_{m-1}$ only
on $X \bigtimes (Z_m \setminus Z_{m-1})$.  
Write $Z_m \setminus Z_{m-1} = \{(p,q)\}$.  
For brevity, we use $g$ in place of $g_{m-1}$.  
For each $x \in X$, we obtain matrices by restricting $g$ as follows:
$$
\align
&A(x) = g |_{\displaystyle \{x\} \bigtimes [p-1,n] \bigtimes [1,q-1]}, \\
&B(x) = g |_{\displaystyle \{x\} \bigtimes \{p\} \bigtimes [1,q-1]}, \\
&C(x) = g |_{\displaystyle \{x\} \bigtimes [p-1,n] \bigtimes \{q\}}.
\endalign
$$
If one of the intervals is empty, the appropriate matrices are zero;
e.g., if $(p,q) = (n,2)$, $C(x) = A(x) = 0$.
\par
For every $x \in X$, we set
$$
\align
K(x) &= B(x)(I-A^*(x)A(x))^{-1/2} \in M_{1, n-p}, \\
L(x) &= (I - A(x)A(x)^*)^{-1/2}C(x) \in M_{q,1}, \\
t(x) &= (I - K(x)K(x)^*)^{-1/2}(I - L(x)^*L(x))^{-1/2}
 \in \Bbb C, \quad t(x) \ge 0, \\
s(x) &= -K(x)A(x)^*L(x) \in \Bbb C.
\endalign
$$
Then, by \cite{DKW, Theorem 1.2}, 
for every number $w$ with $|w| \le t(x)$, the matrix
$$
\pmatrix
B(x) & w+s(x) \\
A(x) & C(x)
\endpmatrix
$$
has norm less than or equal to 1.  In fact, this is also
a necessary condition.  We now define
$$
g_m(x,(i,j)) =
\left\{
\aligned
& g_{m-1}(x,(i,j)), \\
& s(x) + h(t(x), g_{m-1}(x,(p,q)) - s(x)),
\endaligned
\right. \quad
\aligned
&\text{if } (i,j) \ne (p,q), \\
&\text{if } (i,j) = (p,q).
\endaligned
$$
For $(x,(i,j)) \in Y$, if $(i,j) \ne (p,q)$, then
$g_m(x,(i,j)) = g_{m-1}(x,(i,j)) = f(x,(i,j))$.  
If $(x,(p,q)) \in Y$, then condition ($*$),  
together with the properties of $h$, implies that
$$h(t(x), g_{m-1}(x,(p,q)) - s(x)) = 
g_{m-1}(x,(p,q)) - s(x) = f(x,(p,q)) - s(x).$$
Hence $g_m(x,(p,q)) = f(x,(p,q))$ in this case.
This shows that $g_m$ satisfies (1).  
To prove (2), we have to prove ($*$) for $g_m$ and the point
$(i_0, j_0) = (p,q)$.  
But this follows from \cite{DKW, Theorem 1.2}.
\qed \enddemo
\par
The following theorem takes place in the context of
a strongly maximal TAF algebra, $A$, with \cstar-envelope,
$B$ and spectral triple $\spectri$.  Elements of $B$ will 
be viewed as continuous functions on $G$ in the usual way
for groupoid \cstar-algebras.  Also $\MM$ will denote the
collection of all ``finite rectangles'' in $G$; i.e., 
$Q \in \MM$ if $Q = \{(x_i, y_j) \: 1 \le i \le n, 1 \le j \le m\}$
for some $x_i$, $y_j$ in some equivalence class.  
For such a $Q \in \MM$, $T[Q]$ is the matrix obtained by 
restricting $T$ to $Q$ and the norm is the usual matrix norm.
\par
\proclaim{Theorem 6.2} If $\UU$ is a closed $A$-module contained
in $B$ with support set $\s$, then, for any $T \in B$,
$$
\dist(T, \UU) = \sup \bigl\{\|T[Q]\| \: Q \in \MM, Q \cap \s = 
\varnothing \bigr\}
$$
\endproclaim 
\demo{Proof} For every $Q \in \MM$ and $T \in B$, $\|T[Q]\| \le 
\|T\|$; hence 
$$\dist(T,\UU) \ge \sup \{ \|T[Q]\| |: Q \in \MM, Q \cap \s = \varnothing \}. 
	$$
For the other direction, the proof is as in \cite{MS2, Theorem 4.1}
because the previous proposition proves the distance
formula for elementary groupoids of type $n$
and the argument of \cite{MS2, Lemma 4.2} still works
even though the collection $\MM$ here is different
from the corresponding collection there.  (The important
property of sets belonging to $\MM$ is that
$T \longrightarrow T[Q]$ is norm reducing.)
\qed \enddemo
\proclaim{Corollary 6.3} For every ideal $\JJ \subseteq A$
and for every $T \in B$,
$$
\dist(T, \JJ) = \sup \{ \dist(T, \II) \: \II \text{ is a meet
irreducible ideal in } A \text{ and } \II \supseteq \JJ \}.
$$
\endproclaim
\demo{Proof} 
Clearly, $\ge$ holds.  
For the reverse inequality, note that for the left hand side we have
$$\dist(T, \JJ) = \sup \{ \|T[Q]\| \: 
	Q \in \MM, Q \cap \s(\JJ) = \varnothing \}$$
while the right hand side equals
$$ \sup\{\|T[Q]\| \: Q \in \MM, Q \cap \s(\II) = \varnothing 
	\text{ with } \II \supseteq \JJ,
	\II \text{ a meet irreducible ideal }\}.$$
Hence, it is enough to show that if $Q \in \MM$ satisfies 
$Q \cap \s(\JJ) = \varnothing$, 
then there is some meet irreducible ideal $\II \supseteq \JJ$
such that $Q \cap \s(\II) = \varnothing$.  
For this, just assume that 
$$ Q = \{(x_i, y_j) \: x_1 \le x_2 \le \dots \le x_n, 
	y_1 \le y_2 \le \dots \le y_m \} \subseteq [u] \bigtimes [u]
$$
and let $I$ be the interval $[x_1, y_m] \subseteq [u]$.
The meet irreducible ideal $\II$ associated with $I$
satisfies $Q \cap \s(\II) = \varnothing$.  
Since $Q \cap \s(\JJ) = \varnothing$ 
we have also $\s(\II) \supseteq \s(\JJ)$.\qed 
\enddemo

\enddocument